\documentclass[12pt]{article}

\usepackage[T2A]{fontenc}
\usepackage[utf8]{inputenc}

\usepackage{amsfonts}
\usepackage{amssymb}
\usepackage{amsmath}
\usepackage{amsthm}

\def \le {\leqslant}
\def \ge {\geqslant}

\newtheorem{theorem}{Theorem}

\begin{document}

\begin{Huge}
 \centerline{\bf On some open problems in}
\centerline{\bf Diophantine approximation}
\end{Huge}
\vskip+0.5cm

\centerline{ by {\bf Nikolay Moshchevitin}\footnote{ Research is supported by RFBR grant No.12-01-00681-a and by the grant of Russian Government, project 11.  
G34.31.0053 and by the grant NSh-2519.2012.}}
\vskip+0.5cm

\hskip+13.0cm {\it In memory of}

\hskip+11.8cm {\it Vladimir Igorevich Arnold}

\hskip+13.0cm (1937 -- 2010)

\vskip+0.5cm

\begin{small}
 {\bf Abstract.}\, We 
discuss several open problems in Diophantine approximation.
Among them there are famous Littlewood's and Zaremba's conjectures as well as some new and not so famous problems.
 
\end{small}
\vskip+0.5cm
  
 In the present paper I shall give a brief surwey on several problems in Diophantine aproximation
which I was interested in and working on.
Of course 
my briev survey is by no means complete.
 The main purpose of this paper is to make 
the problems under consideration
more popular among  a wider mathematical audience.
 
Basic facts concerning Diophantine approximation one can find in wonderful books
\cite{kok,cassil,ssss2}. As for the Geomerty of Numbers we refer to \cite{cassil2,gl}.
Some topics related to the problems discussed below are considered in recent surveys \cite{Moshew,wald}.

\section{Littlewood conjecture and related problems}
Every paper in Diohantine approximations should begin with the formulation of the Dirichlet theorem  whitch states that
for real numbers $\theta_1,...,\theta_n, n \ge 1$ there exist infinitely many positive integers $q$ such that
$$
\max_{1\le j \le n }||q\theta_j||\le \frac{1}{q^{1/n}}
$$
(here  and in the sequel $||\cdot ||$ stands for
the distance to the nearest integer), or
$$
\liminf_{q\to +\infty} q^{1/n} \,\max_{1\le j \le n}\, ||q\theta_j ||  \le 1
.$$
The famous Littlewood conjecture 
 in Diophantine approximatios supposes that for any two real numbers $\theta_1,\theta_2$ one has
\begin{equation}\label{littel}
\liminf_{q\to +\infty} q\, ||q\theta_1 ||\,||q\theta_2 || =0.
\end{equation}
 The similar multidimensional problem is for given $n \ge 2$ to prove that for any reals $\theta_1,...,\theta_n$ one has
\begin{equation}\label{littel1}
\liminf_{q\to +\infty} q\, ||q\theta_1 || \cdots ||q\theta_n || =0
.\end{equation}
This problem is not solved for any $n \ge 2$. Obviously the statement is false for $n=1$: one can consider 
a badly approximable number $\theta_1$ such that
\begin{equation}\label{littel2}
\inf_{q\in \mathbb{Z}_+} q || q\theta_1|| >0
.\end{equation}
Any quadratic irrrationality satisfies (\ref{littel2});
moreover as it was proved by V. Jarn\'{\i}k \cite{JJ} 
the set of all $\theta_1$ satisfying (\ref{littel2}) has zero Lebesgue measure but  full Hausdorff dimension in $\mathbb{R}$.

Of course Littlewood conjecture is true for almost all pairs $(\theta_1,\theta_2) \in \mathbb{R}^2$.
Moreover, from Gallagher's  \cite{Gal} theorem we know that 
for a positive-valued decreasing to zero function $\psi  (q)$
the inequality
$$
||q\theta_1||\,||q\theta_2||\le \psi (q)
$$
for almost all pairs $(\theta_1,\theta_2)$ in the sense of Lebesgue measure
has infinitely many  solutions in integers $q$ (respectively, finitely many solutions) if the series
$
\sum_q \psi(q) \log q$
diverges (respectively, converges). Thus for almost all $(\theta_1,\theta_2)$ we have
$$
\liminf_{q\to +\infty} q\log^2 q ||q\theta_1||\,||q\theta_2|| =0.
$$

Einsiedler, Katok and Lindenstrauss \cite{EKL} proved that 
the set of pairs $(\theta_1, \theta_2)$ for which (\ref{littel}) is not true is a set of zero Hausdorff dimension
(see also a paper by Venkatesh \cite{veee} devoted to this result).
In my opinion Littlewood conjecture is one of the most exiting open problems in Diophantine approximations.
Some argument for Littlewood conjecture to be true are  given recently by Tao \cite{tao}.

At the beginning of our discussion we would like to formulate Peck's 
theorem  \cite{pekus} concerning approximations to algebraic numbers.

\begin{theorem}\label{pekus}
 Suppose that $n \ge 2$ and $1, \theta_1,...,\theta_n$
form a basis of a real algebraic field of degree $n+1$.
Then there exists a positive constant $C = C(\theta_1,...,\theta_n)$
such that
 there exist infinitely many $q\in \mathbb{Z}_+$ such that simultaneously
$$
\max_{1\le j \le n} ||q\theta_j ||\le
\frac{C}{q^{1/n}}
$$
and 
$$
\max_{1\le j \le n-1} ||q\theta_j ||\le
\frac{C}{q^{1/n}(\log q)^{1/(n-1)}}
.$$

\end{theorem}

Peck's theorem ia a quantitative generalization of a famous theorem by Cassels and Swinnwrton-Dyer \cite{Swin}.
We see that  for a basis of a real algebraic field one has 
$$
\liminf_{q\to +\infty} q \log q\, ||q\theta_1 || \cdots ||q\theta_n || < +\infty,
$$
and so  (\ref{littel1}) is true for these numbers. In particular, 
Littlewood conjecture (\ref{littel}) is  true for numbers $\theta_1,\theta_2$ which form togwther with 
1 a basis of a real cubic field.

We should note here that the numbers $\theta_1,...,\theta_n$ which together with $1$ 
form a basis of an algebraic  field
are simultaneously badly approximable, that is
\begin{equation}\label{bee}
\inf_{q\in \mathbb{Z}_+} q^{1/n}\, \max_{1\le j \le n} ||q\theta_j || >0
\end{equation}
(see \cite{cassil},  Ch. V, \S 3).

A good inrtoduction to Littlewood conjecture one can find in \cite{que}.
Interesting discussion is in  \cite{bullm}.

\subsection{Lattices with positive minima}

Suppose that $1, \theta_1,...,\theta_n$ form a basis of a totally real algebraic field
$\mathbb{K} = \mathbb{Q} (\theta )$ of degree $n+1$. This means that all algebraic conjugates 
$\theta = \theta^{(1)}, \theta^{(2)}, ...,\theta^{(n+1)}$
to $\theta$ are real algebraic numbers.
So 
there exists a polynomial $g_j (\cdot )$  with rational coefficients of degree $\le n$ such that 
$\theta_j  = g_j (\theta)$.
We consider conjugates
$ \theta_j^{(i)} = g_j (\theta^{(i)})$  
and the matrix
$$
\Omega 
=
\left(
\begin{array}{cccc}
1 & \theta_{1}^{(1)}&\cdots&\theta_{n}^{(1)}\cr
\cdots &\cdots &\cdots \cr
1 &\theta_{1}^{(n)}&\cdots&\theta_{n}^{(n)}
\end{array}
\right).
$$
Let
$
G 
$
be a diagonal matrix of dimension $(n+1)\times (n+1)$ with non-zero diagonal elements.
We consider  a lattice of the form 
$$
\Lambda  =  G\,\Omega \,\mathbb{Z}^{n+1}.
$$
Lattices of such a type are known 
as algebraic lattices.

For an arbitrary  lattice $\Lambda \subset \mathbb{R}^{n+1}$
we consider its homogeneous minima
$$
{\cal N} (\Lambda) =
\inf_{ {\bf z} = (z_0,z_1,...,z_{n}) \in \Lambda \setminus \{{\bf 0}\}}\,\,
|z_0z_1\cdots z_{n}|.
$$

One can easiily see that if $\Lambda$ is an algebraic lattice then $ {\cal N}(\Lambda)>0$.

If $ n = 1$ and $\xi, \eta$ are arbitrary badly appoximable numbers (that is satisfying (\ref{bee}))
then for
$$
\Xi
=
\left(
\begin{array}{cc}
1 & \xi\cr
1 &\eta
\end{array}
\right)
$$
the lattice $ L = \Xi \mathbb{Z}^2 \subset \mathbb{R}^2$ will satisfy the property
${\cal N}(L) >0$. Of course one can take $\xi,\eta$ in such a way that $L$ is not an algebraic lattice. So
in the dimension $ n+1 = 2$ there exists a lattice $L\subset \mathbb{R}^2$ which is not an algebraic one, but  ${\cal N}(L) >0$.

A famous {\it  Oppenheim conjecture} supposes that  in the case $n\ge 2$
any  lattice $\Lambda$ satisfying ${\cal N}(L) >0$  is an algebraic lattice.
This conjecture is still open. 
Cassels and Swinnerton-Dyer \cite{Swin} proved that from Oppenheim conjecture in dimension $ n = 2$ Littlewood conjecture 
(\ref{littel}) follows.

Oppenheim conjecture can be reformulated in terms of {\it sails} of lattices (see \cite{herma1,herma2,herma3}).
A sail of  a lattice is a very interesting geometric object whish generalizes Klein's geometric interpretation 
of the ordinary continued fractions algorithm. A 
connection between sails and Oppenheim conjecture was found by Skubenko \cite{sku1,sku2}.
(However the main result of the papers \cite{sku1,sku2} is incorrect: Skubenko claimed the solution of Littlewood conjecture,
howewer he  had a mistake in Fundamental Lemma IV in \cite{sku1}.)

Oppenheim conjecture can be reformulated in terms of closure of orbits of lattices (see \cite{sku1}) and in terms of behaviour of trajectories
of certain dynamical systems. There is a lot of literature related to Littlewood-like problems in lattice theory and dynamical approach 
(see \cite{EKL,kleee,keel,lew,marco,u1,SHA,u2,veee}). 
Some open problems in dynamics related to Diophantie approximations are discussed in \cite{gorodnik}.

\subsection{W.M. Schmidt's conjecture and Badziahin-Pollington-Velani theorem}

For $\alpha, \beta \in [0,1] $ under the condition $ \alpha +\beta = 1 $ and $ \delta >0$
 we consider the sets
$$
{\rm BAD}(\alpha, \beta ;\delta ) = \left\{\xi = (\theta_1,\theta_2 ) 
\in [0,1]^2:\,\,\,\inf_{p\in \mathbb{N}} \,\, \max \{ p^\alpha ||p\theta_1||,
 p^\beta ||p\theta_2||\} \ge \delta \right\}
$$
and
$$
{\rm BAD}(\alpha, \beta  ) = \bigcup_{\delta > 0} {\rm BAD}(\alpha, \beta ;\delta )
 .
$$
  In \cite{SCH1} Schmidt conjectured that for any $\alpha_1,\alpha_2, \beta_1,\beta_2 \in [0,1],\, \alpha_1 +\beta_1 =\alpha_2+\beta_2 =1$
  the intersection
  $$
{\rm BAD}(\alpha_1, \beta_1  ) \bigcap {\rm BAD}(\alpha_2, \beta_2  )
$$
is not empty.
Obviously if Schmidt's conjecture be wrong then Littlewood conjecture be true. But
 Schmidt's conjecture was recently proved in a breakthrough paper by   Badziahin,  Pollington and  Velani \cite{BADZ}.
They proved a more general result: 
\begin{theorem}\label{bpv}
{For any finite collection of  pairs $(\alpha_j,\beta_j),\,\,
0\le \alpha_j,\beta_j \le 1,\,\, \alpha_j+\beta_j = 1,\,\, 1\le j \le r$ 
and for any $\theta_1$ under the condition
\begin{equation}\label{vsop}
\inf_{q\in \mathbb{Z}_+} q||q\theta_1|| >0
\end{equation}
the intersection\begin{equation}\label{jjj}
\bigcap_{j=1}^r \{ \theta_2 \in [0,1]:\,\,\,
(\theta_1, \theta_2 ) \in
 {\rm BAD}(\alpha_j, \beta_j  )
\}
\end{equation}
has full Hausdorff dimension.}
\end{theorem}

Moreover one can take a certain infinite intersection in (\ref{jjj}).

This result was obtained by an original method invented by 
Badziahin, Pollington and Velani.
Author's preprint \cite{mooo} is devoted to an exposition of the method in the simplest case.
The only purpose of the paper \cite{mooo} was to explain the mechanism of the 
method invented by Badziahin, Pollington and Velani.
In this paper 
 it is shown that for
$0<\delta \le 2^{-1622}
$
and $ \theta_1$ such that
 $$
\inf_{q\in \mathbb{N}} q^2||q\theta_1|| \ge\delta,
$$
 there exists $\theta_2$ such that for all integers $A,B$ with $\max (|A|,|B|) >0$  one has
$$
 ||A\theta_1 -B\theta_2 ||\cdot \max(A^2,B^2) \ge {\delta} 
$$
and hence
$$
\inf_{q\in \mathbb{Z}_+}
q^{1/2} \max_{1\le j \le 2 } ||q\theta_j||>0.
$$
This is a quantitative version of a corresponding results from 
\cite{BADZ}.

However the method works with two-dimensional sets only.
Of course we can consider multidimensional BAD-sets.
For example
for $\alpha, \beta, \gamma  \in [0,1] $ under the condition $ \alpha +\beta +\gamma= 1 $ and $ \delta >0$
 one may consider the sets
$$
{\rm BAD}(\alpha, \beta, \gamma ;\delta ) = \left\{\xi = (\theta_1,\theta_2,\theta_3 ) \in [0,1]^3:\,\,\,
\inf_{p\in \mathbb{N}} \,\, \max \{ p^\alpha ||p\theta_1||,\,\,
 p^\beta ||p\theta_2||, \,\, p^\beta ||p\theta_3||\} \ge \delta \right\}
$$
and
$$
{\rm BAD}(\alpha, \beta, \gamma  ) = \bigcup_{\delta > 0} {\rm BAD}(\alpha, \beta , \gamma;\delta )
 .
$$
The question if for given $(\alpha_1,\beta_1, \gamma_1) \neq  (\alpha_1,\beta_1, \gamma_1)$
one has
  $$
{\rm BAD}(\alpha_1, \beta_1 ,\gamma_1 ) \bigcap {\rm BAD}(\alpha_2, \beta_2,\gamma_2  )
\neq \varnothing$$
remains open.

A positive answer is obtained in a very special cases (see, \cite{aaa1,aaa3,aaa2}).

Recently  Badziahin \cite{BADZ1} 
adopted the construction from \cite{BADZ}  to Littlewood-like setting and 
proved the following result.

\begin{theorem}\label{baaad}
The set
$$
\{(\theta_1,\theta_2 )\in \mathbb{R}^2:\,\,\, \inf_{x\in \mathbb{Z}, x \ge 3}
x\log x \,\log\log x  \,||\theta_1 x|| \,
||\theta_2 x|| >0\}
$$
has Hausdorff dimension equal to 2.

Moreover 
if $\theta_1 $ is a badly approximable number
(that is (\ref{vsop}) is valid)
 then the set
$$
\{\theta_2 \in \mathbb{R}:\,\,\, \inf_{x\in \mathbb{Z}, x \ge 3}
x\log x \,\log\log x  \,||\theta_1 x|| \,
||\theta_2 x|| >0\}
$$
has Hausdorff dimension equal to 1.
\end{theorem}

\subsection{$p$-adic version of Littlewood conjecture}

Let for a prime $p$ consider the p-adic norm $|\cdot |_p$, that is
if $  n = p^\nu n_1, (n_1,p)  = 1, \nu\in \mathbb{Z}_+$ then
$|n|_p = p^{-\nu}$. De Mathan and Teuli\'{e} conjectured \cite{demaTH} that for any $\theta \in \mathbb{R}$ one has
$$
\liminf_{q\to +\infty} q|q|_p ||q\theta|| = 0
.
$$
This conjecture is known as {\it p-adic } or {\it mixed} Littlewood conjecture.
Of course the conjecture is true for almost all numbers $\Theta$.
The conjecture can be reformulated as follows: to prove or to disprove that for irrational $\theta$ one has
$$
\inf_{n,q\in \mathbb{Z}_+ }q||p^n q\theta || = 0.
$$
It worth noting that de Mathan and Teuli\'{e} themselves proved \cite{demaTH} that their conjecture is valid for  
 every quadratic irrational
$\theta$.

As  it is shown  in \cite{keel}
from Furstenberg's result 
discussed behind in 4), Subsection 1.4 
it follows that  for distinct $p_1,p_2$
one has
$$
\inf_{m,n,q\in \mathbb{Z}_+ }q||p_1^np_2^m q\theta || = 0
$$
and
so
\begin{equation}\label{newe}
 \liminf_{q\to +\infty} q|q|_{p_1}|q|_{p_2} ||q\theta|| = 0
\end{equation}
for all $\theta$.
Moreover from Bourgain-Lindenstrauss-Michel-Venkatesh's result \cite{JBL}
it follows that for some positive $\kappa$ one has
$$
\liminf_{q\to +\infty} q (\log\log\log q)^\kappa |q|_{p_1}|q|_{p_2} ||q\theta|| = 0.
$$
We should note that the inequality (\ref{newe}) is not the main result of the paper 
\cite{keel} by Einsiedler and Kleinbock.
The main result from \cite{keel} establishes the
zero Hausdorff dimension of the exceptional set in the problem under consideration.

Many interesting metric results and multidimensional conjectures are discussed in \cite{buhh}.

Badziahin and Velani \cite{mixedbad} generalized proved an analog of Theorem \ref{bpv} for mixed Littlewood conjecture:

\begin{theorem}\label{BVV}
The set of reals $\theta$
 satisfying
$$
\liminf_{q\to +\infty} q\,\log q\,\log\log q \,|q|_p\,
||q\theta||>0
$$
has Hausdorff dimension equal to one.
\end{theorem}

In this subsection we consider powers of primes $p^n$ only.
Instead of powers of a prime it is possible to consider other sequences of integers.
This leads to various generalizations. Many interesting results and conjectures of such a kind are discussed in
\cite{bDrmota,meeeba} and \cite{harrapH}.

\subsection{Inhomogeneous problems}

Shapira 
 \cite{SHA}
proved recently two important theorems. We put them below.

\begin{theorem} 
\label{pppp}
Almost all (in the sense of Lebesgue measure) pairs
$(\theta_1,\theta_2) \in \mathbb{R}^2$ satisfy the following property: for every pair
$(\eta_1,\eta_2) \in \mathbb{R}^2$ one has
$$
\liminf_{q\to \infty} q\,\,||q\theta_1-\eta_1||\,\,
||q\theta_2 - \eta_2 || = 0.
$$
\end{theorem}

\begin{theorem}
 The conclusion  of Theorem \ref{pppp}
 is true for numbers $\theta_1, \theta_2$ which form together with $1$
a basis of a totally real algebraic field of degree $3$.
\end{theorem}

Here we should note that the third theorem from \cite{SHA} follows from Khintchine's result (see \cite{HINS}) immediately:

\begin{theorem}\label{rrrr}
 Suppose that  reals $\theta_1$  and $\theta_2$ are
linearly dependent over $\mathbb Z$ together with $1$. Then there
exist reals $\eta_1, \eta_2$ such that
$$
\inf_{x\in \mathbb{Z}_+}\,  x\cdot ||x\theta_1 -\eta_1 || \cdot
||x\theta_2 -\eta_2 ||
>0.
$$
\end{theorem} 
Theorem \ref{rrrr} is discussed in author's paper \cite{inhom}.
Moreover in this paper the author deduces from Khintchine's argument \cite{HINS} the following

\begin{theorem}
 \label{uuuu}
Let $\psi(t) $ be a function increasing to infinity  as $t\to+\infty$.
Suppose that for any $w\ge 1$ we have the inequality
$$
\sup_{x\ge 1}\frac{\psi (wx)}{\psi (x)} <+\infty.
$$
 Then there exist real numbers
$\theta_1, \theta_2$ linearly independent over $\mathbb{Z}$
together with $1$ and real numbers $\eta_1,\eta_2$ such that
$$
\inf_{x\in \mathbb{Z}_+}  x \psi (x)\cdot  ||\alpha_1 x - \eta_1||
\cdot  ||\alpha_2 x - \eta_2|| >0.
$$
\end{theorem}

The following problem is an open one: {\it is it possible that for a constant  function $\psi(t) =\psi_0 >0 \,\,\forall t$
the conclusion of Theorem \ref{uuuu} remains true}.
If not, it means that a stronger inhomogeneous version of Littlewood conjecture is valid.

We would like to formulate here one open problem in imhomogeneous approximations due to Harrap \cite{harrap}.
Harrap \cite{harrap} proved that given $\alpha, \beta \in (0,1),\,\,\alpha + \beta =1$
for a fixed vector
\begin{equation}\label{conde}
(\theta_1,\theta_2) \in {\rm BAD} (\alpha, \beta)
\end{equation}
the set
\begin{equation}\label{seet}
\{ (\eta_1,\eta_2)\in \mathbb{R}^2:\,\,
\inf_q
\max ( q^\alpha ||q\theta_1-\eta_1||, q^\beta||q\theta_2 - \eta_2||>0\}
\end{equation}
 has full Hausdorff dimension.
It is possible to prove that this set is an $1/2$-winning set in $\mathbb{R}^2$
(we discuss winning properties in Subsection 1.6 below).
The following question formulated by Harrap \cite{harrap}:
{\it to prove that the set (\ref{seet}) is a set of full Hausdorff dimension (and even a winning set) 
without the condition (\ref{conde}).}
Of course  in the case $\alpha  =\beta = 1/2$ the positive answer follows from Khintchine's approach
(see results from the paper \cite{mmjm} and the historical discussion there).
However in the case $(\alpha,\beta) \neq (1/2.1/2)$  Harrap's question is still open
\footnote{ A sketch of a proof for Harrap's conjecture is given
in a recent preprint \cite{moshehar}.}.

\subsection{Peres-Schlag's method}

In \cite{PS} Peres and Schlag proved the following result.

\begin{theorem}\label{ps}
 Consider a  sequence $t_n\in \mathbb{R},\, n =1,2,3,..$ Suppose  
that for some $M \ge 2$ one has 
\begin{equation}\label{lus}
\frac{t_{j+1}}{t_j} \ge 1+\frac{1}{M},\,\,\, \forall j \in \mathbb{Z}_+.  \end{equation} 
Then with a certain absolute constant 
  $\gamma
> 0$
for any sequence $\{t_j\}$ under the condition (\ref{lus}) there exists real  $\alpha $ such that
\begin{equation}
||\alpha t_j  || \ge \frac{\gamma}{M\log M},\,\,\, \forall j \in \mathbb{Z}_+.
\label{ququ}
\end{equation}
\end{theorem}
This result has an interesting history with starts from famous Khintchine's paper \cite{HINS}.
We do not want to go into details about this history and refer to pepers \cite{Moshew,MU}.

As it was noted by Dubickas \cite{Db},
an inhomogeneous version of Theorem \ref{ps} is valid: {\it 
with a certain absolute constant 
  $\gamma'
> 0$
for any sequence $\{t_j\}$ under the condition (\ref{lus}) 
and for any sequence of real numbers $\{\eta_j\}$
there exists real  $\alpha $ such that
\begin{equation}\label{peresi}
||\alpha t_j  -\eta_j|| \ge \frac{\gamma}{M\log M},\,\,\, \forall j \in \mathbb{Z}_+.
\end{equation}
 }

One can easily see that for any large integer $M$ there exists an infinite sequence $\{t_j\}$ such that
such that
for any real $\alpha$
there exists infinitely many $j$  with  
\begin{equation}
||\alpha t_j  || \le\frac{1}{M}
\label{ququ1}
\end{equation}
(one may start this sequence $\{t_j\}$ with a finite part $1,2,3,...,M$ and then continue by $1,M, 2M, 3M,..., M\cdot M$, e.t.c.).
Of cource constant $1$ in the numerator of the right hand side may be improved.
The  open question is to find the right {\it order } of approximation in the {\it homogeneous} version of this problem.
In general Peres-Schlag's method does not give optimal bounds. The conjecture is that Theorem \ref{ps}  may be improved on, and 
the optimal result should be stronger than the inequality (\ref{ququ}). 

There are several results and papers dealing with Peres-Schlag's construction
(see  \cite{bumo,bumo1,Db,moNT,MBAD,inhom,MU,PS,RO}).
Here we refer to  four  such results.  

The results in 1) and 2)  below are taken from \cite{inhom}. The paper \cite{inhom} contains some other results related 
to Peres-Schlag's method.

\vskip+0.7cm

1) In  Littlewood-like setting we got the following result.
 
{\it Let $\eta_q, q =1,2,3,..$ be a  sequence
of reals.
Given positive $\varepsilon \le 2^{-14}$ and a badly
approximable real $\theta_1$ such that
$$
||\theta_1 q || \ge \frac{1}{\gamma q}\,\,\, \forall q \in
\mathbb{Z}_+, \,\,\, \gamma >1,$$ there exist $X_0 = X_0
(\varepsilon ,\gamma )$ and a real $\theta_2$ such that
$$
\inf_{q\ge X_0}\,  q \ln^2q \cdot ||q\theta_1|| \cdot ||q\theta_2
-\eta_q||\, \ge \varepsilon.
$$}
If  one consider the sequence $\eta_q =0$ the result behind 
was obtained in \cite{bumo}. It
is worse than Badziahin's Theorem \ref{baaad}.

2) In Schmidt-like settind we proved the following statement.

{\it Suppose that  $\alpha,\beta > 0$ satisfy $\alpha+\beta = 1$. 
 Let $\eta_q, q =1,2,3,..$ be a sequence
of reals.
Let $\eta$ be an arbitrary real number.
 Let $\gamma >0$. Suppose that  $\varepsilon
$ is small enough. Suppose that for a certain   real
$\theta_1$ and for   $ q\ge X_1$  one has
$$
||\theta_1 q|| \ge \frac{\gamma \,(\ln q)^\alpha}{ q^{1/\alpha}}.$$   Then there exist
$X_0 = X_0 (\varepsilon, \gamma  ,X_1 )$ and a real $\theta_2$ such
that
$$
\inf_{q\ge X_0}\,  \max  ((q \ln q)^\alpha \cdot ||q\theta_1-\eta||, \,\,(q\ln q)^\beta \cdot ||q\theta_2
-\eta_q|| )\, \ge \varepsilon.
$$}
Note that
if we take $\eta=0,\eta_q \equiv 0$ we get a result from \cite{MBAD} which is much worse that
Badziahin-Pollilgton-Velani's Theorem \ref{bpv}.

3) For a real $\theta $ we  deal with the sequence $||q^2\theta||$.
Peres-Schlag's argument gives the following statement  (see \cite{moNT})
which solves the simplest problem due to Schmidt \cite{Spoly}.

{\it Given a sequence $\{\eta_q\}$ there exists $\theta$ such that
for all positive integer  $q$ one has
\begin{equation}\label{square}
||q^2\theta- \eta_q||\ge \frac{\gamma}{q\log q}
\end{equation}
(here $\gamma$ is a positive absolute constant).}

Zaharescu \cite{ZA} proved that for any positive $\varepsilon$ and irrational $\alpha $ one has
\begin{equation}\label{square1}
\liminf_{q\to +\infty} q^{2/3-\varepsilon}||q^2\theta|| =0.
\end{equation}
I do not know if this result may be generalized for the value
$$
\liminf_{q\to +\infty} q^{2/3-\varepsilon}||q^2\theta- \eta||
$$
with a real $\eta$.

Nevertheless even in the homogeneous case the lower bound (\ref{square}) is the best known.
So in the homogeneous case we have a gap between (\ref{square}) and (\ref{square1}).

4) F\"urstenberg's sequence.
Consider integers of the form $ 2^m3^m, m, n \in \mathbb{Z}_+$ 
written in the increasing order:
$$
s_0 =1< s_1 =2< s_2=3< s_3 = 4< s_5 = 6< s_6 = 8< s_7 = 9 < s_8 = 12 <...
$$
F\"urstenberg \cite{FU} (simple proof is given in \cite{BOSH}) proved that the sequence of fractional parts
$\{ s_q\theta\}, q\in \mathbb{Z}_+$ is dense for any irrational $\theta$.
Hence for any $\eta$ one has
$$
\liminf_{q\to +\infty}
||s_q\theta -\eta|| =0.
$$
Bourgain, Lindenstrauss, Michel and Venkatesh \cite{JBL} proved a quantitative version of this result.
In particular they show  that if $\theta$
saitisfies for some positive $\beta$ the condition
$$
\inf_{q\in \mathbb{Z}_+}
q^\beta 
||q\theta|| >0
$$
then with some positive $\kappa$
and for any $\eta$
one has
\begin{equation}\label{boura}
\liminf_{q\to +\infty }
(\log\log \log q)^\kappa ||s_q\theta - \gamma || = 0.
\end{equation}
Of course F\"urstenberg had a more general result: instead of the sequence $\{s_q\}$ he considered an arbitrary non-lacunary
multiplicative semigroup in $\mathbb{Z}_+$. The result (\ref{boura}) deals with this general setting also.

Peres-Schlag's method gives the following result:
{\it for an arbitrary sequence $\eta_q, q =1,2,3,...$ there exists irrational $\theta$ such that
$$
\inf_{q\ge 2} \sqrt{q} \log q \,
||s_q\theta -\eta_q || >0
.
$$}
\vskip+0.7cm
One can see that  the results from 1), 2)
with $\eta_q \equiv 0$ are known to be not optimal.
We do not know if the original Theorem  \ref{ps} and the results from 3), 4) with  $\eta_q \equiv 0$ are not optimal.
However I am sure that in {\it  homogemeous} setting the original Theorem  \ref{ps} and the results from 3), 4)
are not optimal and may be improved on.
From the other hand it may happen that the order of approximation in the setting with arbitrary sequence $\{\eta_q\}$
is optimal for some (and even for {\it all}) results  from 1), 2), 3), 4)  and for lacunary sequences.

Of course Peres-Schlag's method 
 gives a thick set 
(a set of full Hausdorff dimension)
of $\theta$'s for which the discussed conclusions hold.

\subsection{Winning sets}

We give the definition of Schmidt's $(\alpha,\beta)$-games and winning sets.
Consider $\alpha,\beta \in (0,1)$,
and a set $S\subseteq \mathbb{R}^d$.
Whites an Blacks are playing the following game.  
Blacks take a closed ball $B_1\subset\mathbb{R}^d$ with diameter $l(B_1)=2\rho$. Then Whites choose a ball
 $W_1\subset B_1$ with diameter $l(W_1)=\alpha l(B_1)$.
Then Blacks choose a ball $B_2\subset W_1$ with diamater
$l(B_2)=\beta l(W_1)$,  and so on... 
In such a way we get a sequence of nested balls
$B_1 \supset W_1 \supset B_2 \supset W_2 \supset \cdots$ with diameters $l(B_i)=2(\alpha\beta)^{i-1}\rho$ and
$l(W_i)=2\alpha(\alpha\beta)^{i-1}\rho$ \,($i=1,2,\dots$).
The set
 $\displaystyle\bigcap_{i=1}^{\infty}B_i=\bigcap_{i=1}^{\infty}W_i$
consists of just one point.
We say that Whites win the game if
 $\smash[t]{\displaystyle\bigcap_{i=1}^{\infty}}B_i\in S$.
A set $S$
is defined to be an
$(\alpha,\beta)$-winning set if Whites can win
the game for any Black's way of playing. 
A set $S$ is defined to be an $\alpha$-winning set  if it is  $(\alpha,\beta)$-winning for every $\beta \in (0,1)$.

 Schmidt \cite{SSSS,ssss1,ssss2}
proved that  for any 
  $ \alpha >0$  an $\alpha$-winning set is a set of full Hausdorff dimension and that the intersection of
a countable family of $\alpha$-winning sets is an $\alpha$-winning set also. 

 For example the set ${\rm BAD}(1/2,1/2)$ is an $1/2$-winning set
(more generally, form Schmidt \cite{ssss1}
we know that the set of badly approximable linear forms is a $1/2$-winning set in any dimension).
Given $\theta\in \mathbb{R}$ the set
$$
\{
\eta \in \mathbb{R}:\,\,
\inf_{q\in \mathbb{Z}_+} q||q\theta - \eta||>0\}
$$
is an $1/2$-winning set \cite{mmjm} (more generally, in \cite{mmjm} there is a result for systems of linear forms).

In 1) - 4) in the previous subsection we discuss the existence of
certain real numbers $\theta_2$ and $\theta$.
In all of these settings it is possible to show that the sets of corresponding $\theta_2$ or $\theta$ have full Hausdorff Dimension.
 Badziahin-Pollington-Velani's Theorem \ref{bpv}, Badziahin-Velani's Theorem \ref{BVV} for mixed Littlewood setting  and  Badziahin's Theorem \ref{baaad}
gives the sets of full Hausdorff dimension also.
However neither in any result from 1) - 4) from the previous subsection, nor in 
 Badziahin-Pollington-Velani's Theorem \ref{bpv}\footnote{Very recently Jinpeng An in  his wonderful paper \cite{JPA}
showed that under the condition $$\inf_q  q^{1/\alpha}||q\theta_1 || >0$$ the set
$$
\{ \theta_2 \in \mathbb{R}:\,\,\,\, (\theta_1,\theta_2) \in {\rm BAD}(\alpha,\beta)\}
$$ is $1/2$-winning.
So he proved the winning property  in Theorem \ref{bpv}. The construction from \cite{JPA} of course gives a better 
quanitiative version of the statement from \cite{mooo} formulated in Section 1.2.

In his next paper \cite{JPA1} Jinpeng An proved that the {\it two-dimensional} set
${\rm BAD} (\alpha,\beta )$ is  $(32\sqrt{2})^{-1}$-winning set.
The construction due to Jinpeng An seems to be very elegant and important.
Probably it can give a solution to  the multi-dimensionsl problem.
As for the constant $(32\sqrt{2})^{-1}$, it seems to me that it may be improved to $1/2$.
}
 and  Badziahin's Theorem \ref{baaad} we do not know if the sets  
constructed are winning.
Moreover we do not know if the set 
${\rm BAD} (\alpha, \beta)$ is a winning set in the case $(\alpha, \beta) \neq (1/2,1/2)$.

The reason is that  the property to be a winning
set is a ``local'' property, but Peres-Schlag's agrument and Badziahin-Pollington-Velani's argument are ``non-local``.
The construction of badly approximable sets by Peres-Schlag and Badziahin-Pollington-Velani
suppose that at a certain level we have a collection of subsegments of a given small segment and we must choose
some good subsegments from the collection. The methods do not 
enable one to say something about the location of good segments form the collection
under consideration. 
All the methods give lower bound for the number of good subsegments in the collection. This is enough to establidh the full Hausdorff dimension,
but does not enough to prove the winning property.

I would be very interesting to study winning properties of the sets arising from Peres-Schlag's method and 
Badziahin-Pollington-Velani's method.

I do not know any   ''natural`` example of a set of badly approximable numbers in a ''natural`` Diophantine problem which has 
full Haudorff dimension but which is not a winning set in the sence of Schmidt games.

At the end of this subsection I would like to fomulate a result by Badziahin, Levesley and Velani \cite{meeeba}:

\begin{theorem}
 \label{jaso1n}
For $\alpha,\beta \in(0,1), \alpha +\beta = 1$ and prime $p$ the set
\begin{equation}\label{odnache}
\{\theta\in \mathbb{R}:\,\,
\inf_{q\in \mathbb{Z}_+}
q \,\max ( |q|_p^{1/\alpha},||q\theta||^{1/\beta}) >0\}
\end{equation}
is $1/4$-winning set.
\end{theorem}

The main result from \cite{meeeba} is more general than Theorem \ref{jaso1n}:
it deals not with p-adic norm $|\cdot|_p$ only but with a norm associated with an arbitrary 
bounded sequence of integers ${\cal D}$.
It is interesting to undestand if the set (\ref{odnache}) is $1/2$-winning.\footnote{recently 
Yaqiao Li \cite{Li} proved that the set  from Theorem \ref{jaso1n} is 1/2-winning indeed.}

\section{ Best approximations}
For positive integers $ m,n $
we consider a real matrix
\begin{equation}\label{ztet} 
\Theta
=\left(
\begin{array}{ccc}
\theta_{1}^1&\cdots&\theta_{1}^{m}\cr
\cdots &\cdots &\cdots \cr
\theta_{n}^{1}&\cdots&\theta_{n}^{m}
\end{array}
\right).
\end{equation}
Suppose that
\begin{equation}\label{good}
\Theta {\bf x} \not\in \mathbb{Z}^n
,\,\,\,\,\, \forall {\bf x}\in \mathbb{Z}^m\setminus \{{\bf 0}\}.
\end{equation}
  We consider a norm $|\cdot |^n_*$ in $\mathbb{R}^n$ and a  norm $|\cdot |^m_*$ in $\mathbb{R}^m$
We are interested mostly in the sup-norm
$$
|{\bf \xi}|^k_{\rm sup} =
\max_{1\le j \le k} |\xi_j|
$$
or in the Euclidean norm
$$
|{\bf \xi}|^k_{2} =
\sqrt{\sum_{j=1}^k |\xi_j|^2}
$$ 
for  a vector  ${\bf \xi} = (\xi_1,...,\xi_k) \in \mathbb{R}^k$
where $k$ is equal to $n$ or $m$.

Let
\begin{equation}\label{bea}
{\bf z}_\nu =({\bf x}_\nu, {\bf y}_\nu)  \in \mathbb{Z}^{m+n},\,\,\,
{\bf x}\in \mathbb{Z}^m,\,\,\, {\bf y} \in \mathbb{Z}^n, \,\,\,\,\,\, \nu =1,2,3,...
\end{equation}
be the infinite sequence of best approximation vectors with respect to the norms $|\cdot |^m_*, |\cdot |^n_*$.
We use the notation
$$ 
M_\nu = |{\bf x}_\nu|^m_*,\,\,\,\,
\zeta_\nu = |\Theta {\bf x}_\nu -{\bf y}_\nu|^n_*.
$$
Recall that the  definition of the best approximation vector can be formulated as follows:
in the set 
$$
\{ {\bf z}_\nu =({\bf x}, {\bf y})\in \mathbb{R}^{m+n}
:\,\,\, | {\bf x}|^m_* \le M_\nu,\,\, |\Theta {\bf x} - {\bf y}|^n_*\le \zeta_\nu \}
$$
there
 is no integer points ${\bf z}_\nu =({\bf x}, {\bf y})$
different from
${\bf 0}, \pm {\bf z}_\nu$.

In this section we formulate some open problems related to the sequence of the best approximations.

\subsection{Exponents of growth for $M_\nu$}

We are interested  in the value
\begin{equation}\label{expo}
G(\Theta) = \liminf_{\nu \to \infty} M_\nu^{1/\nu}.
\end{equation}

First of all we recall well-known general lower bounds.
In \cite{BuLau} it is shown that for sup-norms in $\mathbb{R}^m$ and $\mathbb{R}^n$
for any matrix $\Theta$  under the condition (\ref{good}) one has
$$
G(\Theta ) \ge  2^{\frac{1}{3^{m+n}-1}}
.
$$
This is  a generalization of Lagarias' bound \cite{Laga} for $m=1$.

In fact
Lagarias' result deals with an arbitrary norm: for any norm $|\cdot|^n_*$
 on $\mathbb{R}^n$ and a vector $\Theta \in \mathbb{R}^n$ that has at least one irrational coordinate, 
the inequality 
$$
M_{\nu+2^{n +1}}\ge 2M_{\nu+1} + M_\nu
$$
is true 
for all $\nu \ge 1$. So $G(\Theta)\ge \phi_{1,n}$ where $\phi_{1,n} $ is the maximal
root of the equation $ t^{2^{n+1}} = 2t+1$.

There is another well known statement which is true in any norm. Given a norm $|\cdot|^n_*$ in $\mathbb{R}^n$, consider
the contact number  
$K=K(|\cdot|^n_*)$
 This number is defined as the maximal number of unit balls with
respect to the norm $|\cdot|^n_*$ without interior common points that can touch another unit ball.
 Consider a
 vector $\Theta \in \mathbb{R}^n$ that has at least one irrational coordinate. Then
the inequality 
$$
M_{\nu+K}\ge M_{\nu+1} + M_\nu
$$
holds
for all $\nu \ge 1$. So $G(\Theta)\ge \phi_{2,n}$ where $\phi_{2,n }$ is the maximal root
of the equation $ t^{K} = t+1$.

The general problem is to find optimal bounds for  the value
\begin{equation}\label{infm}
\inf_\Theta G(\Theta )
\end{equation}
for fixed dimensions $m,n$ and for fixed norms  
$|\cdot |^m_*, |\cdot |^n_*$
(the infimum here is taken over matrices $\Theta$ satisfying (\ref{good})).
This problem seems to be difficult.
The only  case where we know the answer is the case $ m=n=1$
(of course in this case there is no dependence on the norms).
For $m=n=1$ the theory of continued fractions gives 
$$ \inf_\Theta G(\Theta ) =G\left(\frac{1+\sqrt{5}}{2}\right) = \frac{1+\sqrt{5}}{2}.$$

Here we consider the case $m=1,n=2$  when  better bounds are known.
In Subsections 2.1.1 and 2.1.2 we formulate the best known results for sup-norm and Euclidean norm.
Any improvement of these bounds may be of interest.
Of course any generalizations to larger values of $n$ are of interest too.
We write
$g_{1,2;\infty}$ for the infimum (\ref{infm}) in the case of sup-norm and   $m=1,n=2$ and
$g_{1,2;2}$ for the infimum (\ref{infm}) in the case of the Euclidean  norm and   $m=1,n=2$ 

\subsubsection{ Case $m=1, n=2$, sup-norm}
In \cite{Moshey} Moshchevitin improved on Lagarias' result form \cite{Laga} by proving 
$$
g_{1,2;\infty}
\ge \phi \cdot\left(\frac{8+13\phi_3}{\phi_3^{13}}
\right)^{\frac{1}{11}} = 1.28040^+,\,\,\, \phi_3 = \sqrt{ \frac{1+\sqrt{5}}{2}}.
$$

\subsubsection{ Case $m=1, n=2$, the Euclidean norm}
Improving on Romanov's result from \cite{Romanov}, Ermakov \cite{Ermakov} proved that
$$
g_{1,2;2} \ge 1.228043.
$$
Here we should note that this result involves numerical computer calculations.

\subsubsection{Brentjes' example related to cubic irrationalities}

Brentjes \cite{Br} considers the following example.
Let $\phi_4 = 1.324^+ $ be the unique real root of the equation
$$
t^3 = t+1.
$$
Consider the lattice
$\Lambda $ consisting of all points of the form
$$
\lambda (\alpha ) =\left(
\begin{array}{c}
\alpha\cr
{\rm Re}\,\alpha '\cr
{\rm Im} \,\alpha '
\end{array}
\right)
,$$
where $\alpha$ is an algebraic integer from the field 
$\mathbb{Q}(\phi_4)$ and $\alpha'$ is one of its algebraic conjugates.
The triple
$$
\lambda (1).\,\,\,
\lambda (\phi_4),\,\,\,\lambda (\phi_4^2)
$$
form a basis of the lattice $\Lambda$.
Brentjes
consideres the sequence of the best approximations
${\bf w}_\nu \in \Lambda , \nu =1,2,3,..$. But his definition differs from our definition behind.
A vector ${\bf w} = (w_0,w_1,w_2)\in \Lambda$ is a best approximation (in Brentjes' sense) if the only points of the lattice
$\Lambda$  belonging to the cylinder
$$
\{ \xi = (\xi_0,\xi_1,\xi_2) \in \mathbb{R}^3:
\,\,|\xi_0|\le |w_0|,\,\,
|\xi_1|^2+|\xi_2|^2
\le
|w_1|^2+|w_2|^2\}
$$
are the points ${\bf 0}, \pm {\bf w}$.
Brentjes shows that these best approximations form a periodic sequence and that
$$
\lim_{\nu \to \infty} |{\bf w}_\nu|^{1/\nu} = \phi_4 = 1,324^+
$$
(here $|\cdot|$ stands for the Euclidean norm in $\mathbb{R}^3$).
This Brentjes' result can be easily obtained by means of the Dirichlet theorem
on algebraic units.

Cusick \cite{couu} studied the  best approximations  for linear form
$$
\phi_4^2x_1+(\phi_4^2-\phi_4)x_2 - y
$$
in the case
$m=2, n=1$ and in sup-norm.

In my opinion the following question remains open:
{\it for the vector $\Theta = \left(
\begin{array}{c}
\phi_4\cr
\phi_4^2
\end{array}
\right)\in \mathbb{R}^2$
(we consider the case $m=1,n=2$)
find the value of 
$G(\Theta)$
defined in (\ref{expo}).
Is it equal to $\phi_4$ or not?}
Probably the solution should be easy.

Lagarias \cite{Laga} conjectured that
in the case $ m=1,n= 2$ for the value $G(\Theta)$ defined in (\ref{ztet}) we have
$$
\inf_{\text{over all norms on }\mathbb{R}^2}\,\,\,
\inf_\Theta G(\Theta) =\phi_4.
$$

\subsection{Degeneracy of dimension: $m=n=2$}

The simplest facts concerning the degeneracy of dimension of subspaces generated by the best approximation vectors are discussed in
\cite{Moshew}.

If
$$
\det \left(
\begin{array}{cc}
\theta_1^1&\theta_1^2\cr
\theta_2^1&\theta_2^2
\end{array}
\right) \neq 0
$$
and (\ref{good}) is satisfied then
for any $\nu_0$ the set of integer vectors $\{{\bf z}_\nu, \nu \ge \nu_0\}$
span the whole space $\mathbb{R}^4$.
So
$$
{\rm dim}\, {\rm span}   \{{\bf z}_\nu, \nu \ge \nu_0\} =4.$$
For a matrix under the condition 
(\ref{good}) the equality
$$
{\rm dim}\, {\rm span}   \{{\bf z}_\nu, \nu \ge \nu_0\} =3$$
never holds.
These facts are proven by Moshchevitin (see Section 2.1 from \cite{Moshew}).

The following question is an opened one.
{\it Does there exist a matrix $\Theta$  (probably, with zero determinant) and satisfying (\ref{good}) such that
for all $\nu_0$ (large enough) one has}
$$
{\rm dim}\, {\rm span}   \{{\bf z}_\nu, \nu \ge \nu_0\} =2\,\, ?$$

\section{Jarn\'{\i}k's  Diophantine exponents}
For a real matrix (\ref{ztet}) satisfying
(\ref{good}) we consider function
$$
\psi_\Theta (t) =
\min_{{\bf x}=(x_1,...,x_m) \in \mathbb{Z}^m: 0< |{\bf x}|^m_{\rm sup}\le t}
\,\,\,
\min_{{\bf y}\in \mathbb{Z}^n}|\Theta {\bf x} -{\bf y}|^n_{\rm sup}
.$$
Sometimes we need to consider the function
$ \psi_{\Theta^*} (t)$ for the transposed matrix $\Theta^*$. In  this case we suppose that $\Theta^*$ satisfies
(\ref{good}) also.
Define ordinary Diophantine exponent
$\omega = \omega(\Theta)$ and uniform Diophantine exponent $\hat{\omega} = \hat{\omega} (\Theta)$:
$$
\omega = \omega(\Theta)=
\sup \left\{
\gamma:\,\,\liminf_{t\to+\infty}
t^\gamma \psi_\Theta (t) 
<+\infty\right\},
$$
$$
\hat{\omega} = \hat{\omega}(\Theta)=
\sup \left\{
\gamma:\,\,\limsup_{t\to+\infty}
t^\gamma \psi_\Theta (t) 
<+\infty\right\}.
$$
In terms of the best approximations (with respect to sup-norm, however here is no dependence on a norm) we have
$$
\omega = \omega(\Theta)=
\sup \left\{
\gamma:\,\,\liminf_{\nu\to+\infty}
M_\nu^\gamma \zeta_\nu
<+\infty\right\},
$$
$$
\hat{\omega} = \hat{\omega(\Theta)}=
\sup \left\{
\gamma:\,\,\limsup_{t\to+\infty}
M_{\nu+1}^\gamma \zeta_\nu 
<+\infty\right\}.
$$
Sometimes we shall use notation $\omega^*(\Theta)$ for $\omega (\Theta^*)$.
There are trivial inequalities
which are valid for all $\Theta$:
$$
\frac{m}{n}\le \hat{\omega}\le \omega \le +\infty.
$$
For $m=1$ one has in addition
$$
\frac{1}{n}\le \hat{\omega}\le 1.
$$
For more details one can see our recent survey \cite{Moshew}.

\subsection{Jarn\'{\i}k's theorems}
In \cite{J} 
 Jarn\'{\i}k proved the following theorem.
\begin{theorem}\label{jahr}
Suppose that $\theta$ satisfies (\ref{good}).

(i) Suppose that $m=1, n\ge 2$ and the column-matrix  $\Theta$ consist at least of two lineqrly independent over $\mathbb{Z}$
together with 1 numbers $\theta_j^1$, Then
\begin{equation}\label{simult}
 \omega \ge  \frac{\hat{\omega}^2}{1-\hat{\omega}}.
\end{equation}

(ii) Suppose that $m=2$. Then
\begin{equation}\label{mtwo}
\omega \ge \hat{\omega}(\hat{\omega} - 1)
.
\end{equation}

(iii) Suppose that $ m \ge 3$ and $ \hat{\omega}\ge (5m^2)^{m-1}$ then
\begin{equation}\label{mthree}
\omega \ge \hat{\omega}^{\frac{m}{m-1}} - 3\hat{\omega} 
.
\end{equation}
 
\end{theorem}

From the other hand
Jarn\'{\i}k proved 
\begin{theorem}\label{jahr1}

(i) Let $m\ge 2$. 
Take real $T >2$.
Then there exists $\Theta$ satisfying (\ref{good}) such that
$$
\omega (\Theta) =  T^m,\,\,\,\, \hat{\omega} (\Theta) = T^{m-1}
.$$

(ii)
Let $m=1, n \ge 2$.
Take real $T>2$  satisfying
$$
T^{n-1} > T^{n-2} + \sum_{k=0}^{n-2} T^k,\,\,\,
T^n > 1 +2\sum_{k=1}^{n-1} T^k.
$$
Then there exists $\Theta$ satisfying (\ref{good}) such that
$$
\hat{\omega} (\Theta) = 
1-\frac{1}{T} -\frac{1}{T^2} - ... - \frac{1}{T^{n-1}}.
$$
$$
{\omega} (\Theta ) = T\,\frac{T^{n-1}-T^{n-2} - ...-1}{T^{n-1}+T^{n-2} + ...+1}
$$
 
\end{theorem}
We see that from (i) of Theorem \ref{jahr1}  it follows that for $\alpha>2^{m-1}$ there exists $\Theta$ such that
$$
\hat{\omega}(\Theta) = \alpha,\,\,\,\omega(\Theta ) = (\hat{\omega}(\Theta))^{\frac{m}{m-1}
}.
$$
From (ii) of Theorem \ref{jahr1} it follows that for $m=1$ and arbitrary $n$ 
for any $\alpha < 1$  close to $1$
there exists a vector $\Theta$ such that
$$
\hat{\omega}(\Theta ) = \alpha,\,\,\,\,
{\omega}(\Theta ) <  \frac{\alpha}{1-\alpha}
.
$$

\subsection{Case $(m,n) = (1,2)$ or $(2,1)$}
Laurent \cite{Lo} proved the following
\begin{theorem}
\label{th47}
The following statemens are valid for the exponents of two-dimensional Diophantine approximations.
 
{(i)} For a vector-row
$\Theta=(\theta^1,\theta^2)\in \mathbb{R}^2$ such that
  $\theta^1,\theta^1$ and 1 are linearly independent over $\mathbb{Z}$ for the values
\begin{equation}
\label{eq92}
w=\hat{\omega}(\Theta),\qquad
w^*=\hat{\omega}(\Theta^*),\qquad
v=\omega(\Theta),\qquad
v^*=\omega(\Theta^*)
\end{equation}
the following statements are valid:
\begin{equation}
\label{eq93}
2\le w\le +\infty,\qquad
w=\frac{1}{1-w^*}\,,\qquad
\frac{v(w-1)}{v+w}\le v^*\le \frac{v-w+1}{w}\,.
\end{equation}

{ (ii)} Given four real numbers $(w,w^*,v,v^*)$,
satisfying~\eqref{eq93} there exists a vector-row
$\Theta=(\theta^1,\theta^2)\in \mathbb{R}^2$,
such that~\eqref{eq92} holds.
\end{theorem}
This theorem is known as ``four exponent theorem''.
It combines together Khintchine's thansference inequalities \cite{HINS} for
ordinary exponents  $\omega, \omega^*$ and 
Jarn\'{\i}k's equality \cite{tiflis} for uniform exponents
$\hat{\omega}, \hat{\omega}^*$ as well as some new results \cite{bulo}.

From Theorem \ref{th47}
it follows that in the case $m=1, n=2$  the inequality (\ref{simult}) is the best possible and cannot be improved.
Also in the case $m=2,n=1$ the inequality (\ref{mtwo}) is the best possible.
The cases  $m=1,n=2$ and $m=2, n=1$ are the only cases when the optimal bounds for $\omega$ in terms 
of $\hat{\omega} $ are known.
In the next two subsections we will formulate  the best known improvements of 
Jarn\'{\i}k's Theorem \ref{jahr}. However all these improvements are far from optimal.
The only possible exception is Theorem \ref{udud} below. The bound of Theorem \ref{udud} may happen to be the optimal one,
however I  am not sure.

To find optimal bounds for $\omega$ in terms of $\hat{\omega} $  (even for specific values of dimensions $m,n$) 
is an interesting open problem.

\subsection{Case $m+n=4$}
Moshchevitin proved the following results
In the case $m=1, n=3$ he get \cite{Czech}

\begin{theorem}\label{udud}
 Suppose that $m=1,n=3$
and the vector $ \Theta=\left( \theta_1,\theta_2, \theta_3\right)$
consists of numbers linearly independent, together with 1, over
$\mathbb{Z}$. Then
\begin{equation}\label{shmu}
\omega  \ge \frac{\hat{\omega}}{2} \left( \frac{\hat{\omega}}{1-\hat{\omega}}
+\sqrt{\left(\frac{\hat{\omega}}{1-\hat{\omega}}\right)^2
+\frac{4\hat{\omega}}{1-\hat{\omega}}}\right) .
\end{equation}.

\end{theorem}

The inequality (\ref{shmu}) is better than Jarn\'{\i}k's inequality (\ref{simult}) for all values of $\hat{\omega}(\Theta)$.

In \cite{Mo1,Moshew} the following two theorems are proved
(a proof of Theorem \ref{uo1} was just sketched).

\begin{theorem}\label{uo}
 Suppose that $m=3,n=1$ and the matrix $\Theta = (\theta^1,\theta^2,\theta^3)$ consists of numbers
linearly independent over   $\mathbb{Z}$  together with $1$. Then
\begin{equation}\label{novoe}
\omega \ge \hat{\omega}  \cdot \left(
\sqrt{\hat{\omega}  +\frac{1}{\hat{\omega}  ^2}-\frac{7}{4}}+\frac{1}{\hat{\omega}  }-\frac{1}{2}\right).
\end{equation}

\end{theorem}

The inequality (\ref{novoe}) is better than Jarn\'{\i}k's inequality (\ref{mthree}) for all values of $\hat{\omega}(\Theta)$.

\begin{theorem}\label{uo1}
 Consider four real numbers $\theta_j^i,\,\, i,j=1,2$ linearly independent over $\mathbb{Z}$ together with 1.
Let $m=n=2$ and consider the matrix
$$
\Theta =
\left(
\begin{array}{cc}
\theta^1_1&\theta^2_1\cr
\theta^1_2&\theta^2_2
\end{array}
\right)
$$
satisfying (\ref{good}).
Then
\begin{equation}
\label{2mod}
\omega \ge
\frac{1-\hat{\omega} +\sqrt{(1-\hat{\omega})^2+4\hat{\omega} (2\hat{\omega}^2-2\hat{\omega}+1)}}{2}
.
\end{equation}
\end{theorem}

The inequality (\ref{2mod})
 improves on the inequality (\ref{mtwo})  for $\hat{\omega} (\Theta) \in \left(1,\left(\frac{1+\sqrt{5}}{2}\right)^2\right)$.\footnote{Recently  the author \cite{dfa} improved Theorem \ref{uo1}.
Now the result is better than (\ref{mtwo})  for all admissible values of $\hat{\omega}$}

It may happen that Theorem \ref{udud}
gives the optimal bound.\footnote{Very recently the author proved that Theorem  \ref{udud}  gives the optimal bound.}
I am sure that the inequality from Theorem \ref{uo1} may be improved.

\subsection{A result by W.M. Schmidt and L. Summerer (2011)}

Very recently Schmidt and Summerer \cite{ss} 
improved on 
Jarn\'{\i}k's Theorem \ref{jahr} in the cases $m=1$ and $n=1$.

For $m = 1$ and arbitrary $n\ge 2$ they obtained the bound
\begin{equation}\label{oop}
\omega  \ge
\frac{\hat{\omega}^2 + (n-2)\hat{\omega} }{(n-1)(1-\hat{\omega})}
.
\end{equation}
As for the dual setting with $n=1$ and $m\ge 2$ they proved that
\begin{equation}\label{oopp}
\omega \ge (m-1) \frac{\hat{\omega}^2-\hat{\omega}}{1+(m-2)\hat{\omega}} 
.
\end{equation}
No analogous inequalities are known in the case when both $n$ and $m$ are greater than one.

The result by Schmidt and Summerer deals with successive minima for 
one-parametric families of lattices and relies on their earlier research
\cite{SSS} and Mahler's theory of compound and pseudocompound bodies \cite{maa}.
We suppose to write a separate paper concerning Schmidt-Summerer's result and its possibble extensions,
jointly with O. German.

Here we should note that in the cases $m=1,n=3$ and $m=3,n=1$ the inequalities (\ref{shmu}) and ({\ref{novoe}) are better than
(\ref{oop}) and (\ref{oopp}), correspondingly.

\subsection{Special matrices}

Here we would like to formulate one open problem which seems to be not too  difficult.
Consider a special set $\hbox{\got W}$  of matrices $\Theta$.
A matrix $\Theta$ belongs to $\hbox{\got W}$ if (\ref{good}) holds and moreover  there exists infinitely many
$(m+n)$-tuples of
{\it consecutive} 
best approximation vectors
$$
{\bf z}_{\nu}, {\bf z}_{\nu+1}, {\bf z}_{\nu+2} ,...,{\bf z}_{\nu + m+n-1}
$$
consisting of {\it  linearly independent } vectors in $\mathbb{R}^{m+n}$.
The definitoion of ${\bf z}_j \in \mathbb{Z}^{m+n} $ (see (\ref{bea})) is given in the very beginning of Section 2.

I think that it is not dificult to improve on  all the inequalities from  
Jarn\'{\i}k's Theorem \ref{jahr} in the case $\Theta \in \hbox{\got W}$.\footnote{In fact,
it can be easily seen that for 
$\Theta \in \hbox{\got W}$ one has
$$
\omega (\Theta) \ge G \cdot \hat{\omega} (\Theta),
$$
where $G$ is the largest root of the equation
$$
-\sum_{j=1}^{n-1}\frac{\hat{\omega}}{x^j}+1-\hat{\omega} + \sum_{j=1}^{m-1}x^j = 0.  
$$
For $m+1,n=1$ the value of $G=\frac{\alpha}{1-\alpha}$ coinsides with the analogous value from (\ref{simult}) from
Jarn\'{\i}k's Theorem
\ref{jahr}. In the case $m=1,n=3$ the value of $G=\frac{1}{2} \left( \frac{\hat{\omega}}{1-\hat{\omega}}
+\sqrt{\left(\frac{\hat{\omega}}{1-\hat{\omega}}\right)^2
+\frac{4\hat{\omega}}{1-\hat{\omega}}}\right)$ coinsides with those from author's Theorem \ref{udud}.}

Moreower I think that in this case it is possible to get optimal inequalities and to prove the optimality of these inequalities by
constructing special matrices $\Theta \in \hbox{\got W}$.

\section{ Positive integers}

In this section we consider collections of  real numbers $\Theta =(\theta^1,...,\theta^m), \, m\ge 2$.
(The index $n=1$ is omitted here.)
We are interested in small values of the linear form
$$
||\theta^1 x_1+...+\theta^nx_n
||$$
in positive integers $x_1,..., x_n$.
Put
$$
\psi_+(t) 
=\psi_{+;\Theta } (t) =\min_{x_1,...,x_m \in \mathbb{Z}_+,\,\, 0< \max (x_1,...,x_m) \le t}
||\theta^1x_1+...+\theta^mx_m||.
$$
We introduce Diophantine exponents
$$
\omega_+ = \omega_{+} (\Theta) =
\sup\{ \gamma:\,\,\,
\liminf_{t \to \infty} t^\gamma \psi_{+;\Theta} (t) < \infty\},
$$ 
and
$$
\hat{\omega}_+ = \hat{\omega}_{+} (\Theta) =
\sup\{ \gamma:\,\,\,
\limsup_{t \to \infty} t^\gamma \psi_{+;\Theta} (t) < \infty\}.
$$

\subsection{The case $m=2$: W.M. Schmidt's theorem and its extensions}

Put
$$
\phi = \frac{1+\sqrt{5}}{2} = 1.618^+.
$$
In  1976 
W.M. Schmidt\cite{SCH} proved the following theorem.
\vskip+0.5cm
\begin{theorem}[W.M. Schmidt]\label{uhu}
Let real numbers $\theta^1,\theta^2$
be linearly independent over $\mathbb{Z}$
together with 1. Then there exists a sequence of integer two-dimensional vectors
 $(x_1(i), x_2(i))$
 such that

 1.\,\, $x_1(i), x_2(i) > 0$;

2.\,\, $||\theta^1x_1(i)+\theta^2x_2(i) ||\cdot (\max \{x_1(i),x_2(i)\})^\phi \to 0$ as $ i\to +\infty$.
\label{shpos}
\end{theorem}
\vskip+0.5cm
In fact W.M. Schmidt proved (see discussion in \cite{BUKR}) that for $ n = 2$ for $\Theta =  (\theta^1,\theta^2)$ under consideration
 one has the inequality
\begin{equation}\label{knipers}
 \omega_{+} \ge \max \left( \frac{\hat{\omega}}{\hat{\omega}-1}; \hat{\omega} - 1 +\frac{\hat{\omega}}{\omega}\right)
\end{equation}
from  which
we deduce
\begin{equation}\label{knipers1}
 \omega_{+} (\Theta)
\ge \phi
 .
\end{equation}
From Schmidt's argument one can easily see that for $\theta^1,\theta^2$ linearly independent together with 1 one has 
\begin{equation}\label{knipers}
 \hat{\omega}_{+} \ge   \frac{{\omega}}{{\omega}-1}.
\end{equation}
We would like to note here that Thurnheer (see Theorem 2 from \cite{aa1989}) showed that for
\begin{equation}\label{knee}
\frac{1}{2}\le \omega^* = \omega^* (\Theta) \le 1 
\end{equation}
($\omega^*(\Theta)$ was defined in the beginning of Section 3, here it is the Diophantine exponent for simultaneous approximations for numbers
$\theta^1,\theta^2$)
one has
\begin{equation}\label{knipers2}
 {\omega}_{+} \ge \frac{\omega^*+1}{4\omega^*} +\sqrt{\left(\frac{\omega^*+1}{4\omega^*} \right)^2+1}.
\end{equation}
(inequality  (\ref{knipers2}) is a particular case  of a general result obtained by Thurnheer).

A lower bound for $\omega_+$ in terms of $\omega  $ was obtained by the author in \cite{moshe}. It was based on the original Schmidt's 
argument from \cite{SCH}. However the choice of parameters in  \cite{moshe} was not optimal. Here we explain
the optimal choice \cite{m2012}.
From Schmidt's proof 
and Jarn\'{\i}k's result (\ref{mtwo})
one can easily see that
$$
\omega_{+} \ge
\max 
\left\{ g:\,\,\,
\max_{y,z \ge 1:\,\, y^{\hat{\omega}-1}\le z\le y^{\omega/\hat{\omega}}}
\,
\max_{y^{-\omega}\le x\le z^{-\hat{\omega}}}\,\,
\min\left(
x^{1-g}z^{-g}; x y^{-1}z^{g+1}
 \right) \le 1 
\right\}
.$$
The right hand side here can be easily calculated.
We divide the set
$$
\hbox{\got A} =
\left\{ (\omega, \hat{\omega}) \in \mathbb{R}^2:\,\,\, 
\hat{\omega}\ge 2,\,\, \omega \ge \hat{\omega}(\hat{\omega} -1) \right\}
$$
of all admissible values of $(\omega,\hat{\omega})$ into two parts:
$$
\hbox{\got A} =
\hbox{\got A}_1\cup \hbox{\got A}_2,
$$ 
$$
\hbox{\got A}_1 =
\left\{ (\omega, \hat{\omega}) \in \mathbb{R}^2:\,\,\,
2\le
\hat{\omega}\le \phi^2,\,\, \omega \ge \frac{\hat{\omega}(\hat{\omega}-1)}{3\hat{\omega}-\hat{\omega}^2 -1}\right\},
$$
$$
\hbox{\got A}_2 = \hbox{\got A} \setminus \hbox{\got A}_1
.$$
If $(\omega,\hat{\omega}) \in \hbox{\got A}_1$
then
$$
\omega_+ \ge G(\omega) = \frac{1}{2}\left(\frac{\omega+1}{\omega}+\sqrt{\left(\frac{\omega+1}{\omega}\right)^2+4}\right)
$$
(the function $ G(\omega)$ on the right hand side decreases from $G(2) = 2$ to $G(+\infty ) = \phi$).
If 
$(\omega,\hat{\omega}) \in \hbox{\got A}_2$ then
\begin{equation}\label{ooa}
\omega_+ \ge \hat{\omega} - 1+\frac{\hat{\omega}}{\omega}
\end{equation}
So
\begin{equation}\label{fff}
\omega_+ \ge
\max\left(\frac{1}{2}\left(\frac{\omega+1}{\omega}+\sqrt{\left(\frac{\omega+1}{\omega}\right)^2+4}\right);
\hat{\omega} - 1+\frac{\hat{\omega}}{\omega}\right)
,\end{equation}
and this is the best bound in terms of $\omega, \hat{\omega}$ which one can deduce from Schmidt's argument from \cite{SCH}.

\subsection{A counterexample to W.M. Schmidt's conjecture
}

In the paper \cite{SCH} W.M. Schmidt wrote that he did not know if the exponent $\phi$ in Theorem \ref{shpos}
may be replaced by a lagrer constant. At that time he was not
able even to rule a possibility
that there exists an infinite sequence $(x_1(i),  x_2(i))\in \mathbb{Z}^2$ with condition 1. and such that
\begin{equation}\label{opip}
||\theta^1x_1(i)+\theta^2x_2(i) ||\cdot (\max \{x_1(i),x_2(i)\})^2 \le c(\Theta)
\end{equation}
with some large positive $c(\Theta)$.
Later in \cite{SCH1} he conjectured that the exponent $\phi$ may be replaced by any exponent of the form
$2 -\varepsilon, \varepsilon >0$ and wrote that probably such a result should be obtained
by analytical tools.
 It happened that this conjecture is not true. 
 In  \cite{Mopositi} the author proved the following result.
\begin{theorem}\label{moshepos1}
Let $\sigma = 1.94696^+$ be the largest real root of the equation
$
 x^4 - 2x^2-4x+1=0.
$ There exist real numbers
 $\theta^1,\theta^2$
such that they are linearly independent over $\mathbb{Z}$
together with 1 and for every  integer vector
$(x_1,x_2)\in\mathbb{Z}^2$ with $x_1 , x_2 \ge 0$  and
$\max (x_1,x_2) \ge 2^{200}$
one has
$$
||\theta^1x_1+\theta^2x_2||\ge \frac{1}{2^{300}(\max (x_1,x_2))^\sigma}.
$$
\end{theorem}
Theorem \ref{moshepos1} shows that W.M. Schmidt's conjectue discussed in previous subsection turned out to be  false.

Here we should note that for the numbers constucted in Theorem \ref{moshepos1} one has
$$
\omega = \frac{(\sigma+1)^2(\sigma^2-1)}{4\sigma} = 3.1103^+ ,
\,\,\,
\hat{\omega} = \frac{(\sigma+1)^2}{2\sigma} = 2.2302^+ .
$$
So $(\omega, \hat{\omega}) \in \hbox{\got A}_2$ and the inequality (\ref{ooa}) gives
$$
\omega_+ \ge \frac{\sigma +2}{\sigma^2-1} = 1.413^+.
$$
However from the proof of Theorem \ref{moshepos1} (see \cite{Mopositi}) it is clear
that for the numbers constructed one has $\omega_+ = \sigma = 1.94696^+$.

\subsection{$m =2$: large domains}
 
The original paper \cite{SCH} contained 5 remarks related to Theorem \ref{shpos}.
One of these remarks was as follows.
The condition 1. in Theorem \ref{shpos}  may be
replaced by
a condition $|\alpha_{1,1,}x_1(i)+\alpha_{1,2}x_2(i)| < |\alpha_{2,1}x_1(i)+\alpha_{2,2} x_2(i)|
$
where $ \alpha_{1,1}\alpha_{2,2}-\alpha_{1,2}\alpha_{2,1}\neq 0$.
In this new setting we deal with good approximations from  an ``angular
domain''.
Later Thurnheer \cite{T}
got a result dealing with even larger domain.
For  positive parameters $\rho, \tau$ he considered the domain
\begin{equation}\label{doma}
 \Phi_0(\rho, \tau)
=
\{ (x_1,x_2) \in \mathbb{R}^2:\,\, |x_2|\le |x_1|^\rho\} \cup
\{ (x_1,x_2) \in \mathbb{R}^2:\,\, |x_1|\le |x_2|^\tau\} 
\end{equation}
and its image $\Phi (\rho, \tau)$ under a non-degenerate linear transform.
Thurnheer \cite{T} proved the following 
\begin{theorem}\label{th1}
Suppose that parameters  $\rho> 1,
\tau \ge 0$ and $1<t\le r\le 2 $ satisfy the condition
\begin{equation}\label{cococo}
 (1-\tau) (\rho(t^2r - tr-t-r-1)+t^2)+(1-\rho)(t^2-1)\le 0.
\end{equation}
Then  there exist infinitely many integer points $(x_1,x_2)$  such that
$$
(x_1, x_2) \in \Phi (\rho, 0)\,\,\,\,\,\,\,\,\text{and}\,\,\,\,\,\,\,\,
||\theta^1x_1+\theta^2x_2||\le c_1 (\max (|x_1|, |x_2|))^{-r},
$$
or
$$
(x_1, x_2) \in \Phi (1, \tau)\,\,\,\,\,\,\,\,\text{and}\,\,\,\,\,\,\,\,
||\theta^1x_1+\theta^2x_2||\le c_2(\max (|x_1|, |x_2|))^{-t}.
$$
Here $c_{1,2}$ are positive constants depending on $\Theta$ and $\rho, \tau$.
\end{theorem}

Thurnheer \cite{T} considered three special cases of his Theorem \ref{th1}:

1. by putting 
\begin{equation}\label{roo1}
 \rho = 7/4, \tau = 0
\end{equation}
and $ t= r = 2$
one can see that there exist infinitely many  integer vectors $(x_1, x_2)$ such that
\begin{equation}\label{qq}
 (x_1,x_2) \in \Phi (7/4,0)\,\,\,\,\,\,\,\,
\text{and}\,\,\,\,\,\,\,\,
||\theta^1x_1+ \theta^2x_2||\le c_3 (\max (|x_1|, |x_2|))^{-2}
\end{equation}
(with a certain value of $c_3>0$); 

2. by putting 
\begin{equation}\label{rooo}
 1<\rho \le 7/4,\,\,\,\,
\tau = \frac{7-4\rho}{4-\rho}
\end{equation}
and $ t= r = 2$
one can see that there exist infinitely many  integer vectors $(x_1, x_2)$ such that
\begin{equation}\label{qqq}
 (x_1,x_2) \in \Phi (\rho,\tau)\,\,\,\,\,\,\,\,
\text{and}\,\,\,\,\,\,\,\,
||\theta^1x_1+ \theta^2x_2|| \le c_4 (\max (|x_1|, |x_2|))^{-2}
\end{equation}
(with a certain value of $c_4>0$);

3. for any $\rho \in  (1,7/4]$ and $\tau = 0$ one can consider the largest root $s(\rho)$ of the equation
\begin{equation}\label{roooq}
 \rho x^3 - 2(\rho -1) x^2 - 2\rho x - 1 =0.
\end{equation}
Then
one can see that there exist infinitely many  integer vectors $(x_1, x_2)$ such that
\begin{equation}\label{roo12}
 (x_1,x_2) \in \Phi (\rho,0)\,\,\,\,\,\,\,\,
\text{and}\,\,\,\,\,\,\,\,
||\theta^1x_1+ \theta^2x_2|| \le c_5 (\max (|x_1|, |x_2|))^{-s(\rho)}
\end{equation}
(with a certain value of $c_5>0$); 
note  that  $s(1) = \phi =\frac{1+\sqrt{5}}{2}$, and this gives Schmidt's bound (\ref{knipers1}).

\subsection{Large dimension ($m >2$)
 }

\subsubsection{A remark related to Davenport-Schmidt's result}

Another remark  to Theorem 1 from \cite{SCH} tells us that  
``no great improvement is affected by allowing a large
number of variables''.
W.M. Schmidt showed the following result to be true.

\begin{theorem}
\label{th56}
There exists a vector $\Theta=(\theta^1,\dots,\theta^m)$, $m\ge 3$ such that:

-- $1, \theta^1,...,\theta^m$ are linearly independent over $\mathbb{Z}$;

-- for any positive~$\varepsilon$ there exists a positive~$c(\varepsilon)$ such that
$$
\|\theta^1x_1+\dots+\theta^m x_m\|>c(\varepsilon)
\Bigl(\,\max_{1\le i\le m}|x_i|\Bigr)^{-2-\varepsilon}
$$
for all integers $x_1,\dots,x_k$ under the condition $x_i >0$, \,$i=1,\dots,k$.
\end{theorem}
To prove Theorem~\ref{th56} one should use a result by H. Davenport and W.M. Schmidt from the paper~\cite{92}.
This result is based on existence of very singular  vectors. 
Theorem \ref{th56} shows that  for $ m\ge 3$ for linearly independent collection 
$\Theta$
it may happen that
\begin{equation}\label{may}
\omega_{+} (\Theta ) \le 2
.
\end{equation}

\subsubsection{ General Thurnheer's lower bounds}

Here we formulate three general results by Thurnheer from \cite{aa1989}.
Its particular case (inequality (\ref{knipers2})) was discussed above.
Thurnheer used the Euclidean norm to formulate his result. Of course it is not of importance and we may use sup-norm.

Given $\varepsilon >0$ consider the domain
$$
\Psi = \Psi_\varepsilon =
\left\{ (x_1,...,x_m) \in \mathbb{R}^m:\,\,
|x_m|\le \varepsilon \max_{1\le j\le m-1}|x_j|\right\}.
$$

\begin{theorem}
 \label{turn}
Suppose that $1,\theta^1,...,\theta^m$ are linearly independent over $\mathbb{Z}$.
Put
$$
v(m) =\frac{1}{2}
\left(
m-1+\sqrt{m^2+2m-3}
\right)> m - \frac{1}{m}.
$$
Then  there exists infinitely many integer vectors
$
(x_1,..,x_m) \in \Psi
$
such that
$$
||\theta^1x_1+...+\theta^mx_m|| \le \delta 
(\max_{1\le j\le m-1}|x_j|)^{-v(m)}
$$
(here $\delta$ is an arbitrary small fixed positive number).
\end{theorem}

Another Thurnheer's result deals with aproximations from a larger domain.
For $w>0$ put
$$
\Phi (w)=
\left\{
(x_1,...,x_m) \in \mathbb{R}^m:\,\,
|x_m|\le (1+\varepsilon )\left(\sum_{j=1}^{m-1}|x_j|^2\right)^{w/2}\right\}\cup$$
$$
\cup\left\{
(x_1,...,x_m) \in \mathbb{R}^m:\,\, \sum_{j=1}^{m-1}|x_j|^2\le 1
\right\},\,\,\, \,\, \varepsilon >0.
$$

\begin{theorem}
 \label{turne}
Put
$$
w=
w(m) =1+\frac{1}{m}+\frac{1}{m^2}.
$$
Then for any real $\Theta$ and
and for any positive $\delta$ there exists infinitely many integer vectors
$
(x_1,..,x_m) \in \Phi (w)
$
such that
$$
||\theta^1x_1+...+\theta^mx_m|| \le (1+\delta ) 
(\max_{1\le j\le m-1}|x_j|)^{-m}
.
$$
\end{theorem}

Another result deals with a lower bound in terms of $\omega^*$.

\begin{theorem}\label{tu}
Suppose that $1,\theta^1,...,\theta^m$ are linearly independent over $\mathbb{Z}$.
Suppose that
\begin{equation}\label{uu}
\frac{1}{m}<
\omega^* = \omega^*(\Theta)\le \frac{1}{m-1}.
\end{equation}
Put
$$
u_0 (m,\omega^*) =
\frac{1}{2m\omega^*}\left(
\omega^* (m-1)^2 +1+ \sqrt{(\omega^* (m-1)^2+1)^2 + 4m^2(m-1)(\omega^*)^2}
\right)
$$
 Then 
for any $u< u_0 (m,\omega^*)$
there exists infinitely many integer vectors
$
(x_1,..,x_m) \in \Psi
$
such that
$$
||\theta^1x_1+...+\theta^mx_m|| \le 
(\max_{1\le j\le m-1}|x_j|)^{-u}
.
$$

\end{theorem}

\subsubsection{From Thurnheer to Bugeaud and Kristensen}
Bugeaud and Kristensen \cite{BUKR} considered the following Diophantine exponents.
Let $ 1\le l\le m$. Consider the set
$$
\Psi (m,l) = \Psi_\varepsilon (m,l)=
\left\{
(x_1,...,x_m \in \mathbb{R}^m:
\,\,\,
\max_{l+1\le j \le m} |x_j|\,\le  \max_{1\le j \le l} |x_j| 
\right\}.
$$
Diophantine exponent  $\mu_{m,l} = \mu_{m,l}(\Theta)$
is defined as the supremum over all $\mu$  such that the inequality
$$
||\theta^1x_1+...+\theta^m x_m||\le (\max_{1\le j \le m}|x_j|)^{-\mu}
$$
has infinitely many solutions in $$ (x_1,...,x_m) \in \Psi(m,l) \cap \mathbb{Z}^m .$$
So Diophantine exponent $\mu_{m,1}$ corresponds just  to $\omega_{+}$.
Thurnheer's Theorems \ref{turn} and \ref{tu}
have the following interpretation in terms of $\mu_{m,m-1}$.

Suppose that $1,\theta^1,...,\theta^m$ are linearly independent over $\mathbb{Z}$. Then
\begin{equation}\label{vem}
\mu_{m,m-1} \ge v (m).
\end{equation}
If in addition  (\ref{uu}) holds then
\begin{equation}\label{vem1}
\mu_{m,m-1}\ge u_0 (m,\omega^*).
\end{equation}

Bugeaud and Kristensen formulate the following result.

\begin{theorem}\label{tubuk}
Suppose that $1,\theta^1,...,\theta^m$ are linearly independent over $\mathbb{Z}$. Then
$$
\mu_{m,l} \ge \frac{l\hat{\omega}}{\hat{\omega} -m+l}
$$
and
$$
\mu_{m,m-1} \ge \hat{\omega} -1+\frac{\hat{\omega}}{\omega}.
$$
\end{theorem}
In fact linearly indenendency condition here is necessary.

Of course from this theorem the bound
(\ref{vem}) follows immediatelly.

Here we should note that the main results of the paper \cite{BUKR} deal with metric prorerties of exponents $\mu_{m,l}$.
Also in \cite{BUKR} several interesting problems are formulated.

\subsection{ Open questions}

1. What is the optimal exponent $\inf_{\theta^1,\theta^2 -\,\text{independent}} \omega_{+}(\Theta)$ in the problem for a linear form in two positive variables? 
Is it $\phi$ or $\sigma$ or something else between 
$\phi$ and $\sigma$?\footnote{ Recently Damien Roy \cite{DROY}
anounced that the exponent $\phi = \frac{1+\sqrt{5}}{2} = 1.618^+$
from Schmidt's Theorem \ref{uhu} is optimal.}

2.
What are the best possible lower bounds for
 $\omega_{+}$ and $\hat{\omega}_{+}$ in terms of $\omega$, $\omega^*$ and  $\hat{\omega}$?
Any improvement of any of the  lower bounds (\ref{knipers}, \ref{knipers2}, \ref{fff}) will be of interest, in my opinion.
Of course any improvement of  lower bounds  for $\mu_{m.l}$ given in (\ref{vem}, \ref{vem1}) as well as of the bounds from  Theorem \ref{tubuk} 
will be of interest.

3. As it was shown behind for $\theta^1,...,\theta^m, m\ge 3$ linearly independent over $\mathbb{Z}$ together with 1 it may happen
(\ref{may}).
But in view of  Theorem \ref{moshepos1}
I may conjecture that for $m = 3$ (or even for an arbitrary $m$) there exist  a collection of linearly independent numbers
$1,\theta^1,...,\theta^m$ such that $ \omega_{+}(\Theta) < 2$.

4. What are optimal exponents in the Thurnheer's setting for large domains $\Phi, \Psi$?
In particular, are the values of parameters $\rho, \tau$ from (\ref{roo1}) and (\ref{rooo}) optimal
to get (\ref{qq}) and (\ref{qqq})
 or not?
What is the optimal values of $s(\rho)$ to conclude that (\ref{roo12}) has infinitely many solutions in integers $(x_1,x_2)$?
Any improvements of the discussed results is of interest.
Similar questions may be formulated for multi-dimensional results.

In view of our Theorem \ref{moshepos1} I think that it is possible  to solve  some of the problems formulated in this subsection.

\section{Zaremba conjecture} 

For an irreducible rational fraction $\frac{a}{q}\in\mathbb{Q}$ we consider its continued fraction expansion 
$$
\frac{a}{q}=[b_0;b_1,\dots,b_s]= 
$$
\begin{equation}\label{exe}
=
b_0  +
\frac{1}{\displaystyle{b_1+\frac{1}{\displaystyle{b_2 +
 \frac{1}{\displaystyle{b_3 +\dots +
\displaystyle{\frac{1}{b_{s}}} }}}}}} ,\,\,\,\,\, b_j =b_j(a) \in \mathbb{Z}_+,\,\, j\ge 1
.
\end{equation}
The famous Zaremba's conjecture \cite{zara} supposes that there exists an absolute constant $\hbox{\got k}$
with the following property:
for any positive integer $q$ there exists $a$ coprime to  $q$ such that in the continued fraction expansion (\ref{exe}) all partial quotients are bounded:
$$ b_j (a) \le \hbox{\got k},\,\, 1\le j  \le s = s(a).
$$
In fact Zaremba conjectured that $ \hbox{\got k}=5$. Probably 
for large prime $q$  even $ \hbox{\got k}=2$ sould be enough, as it was conjectured by Hensley .

\subsection{What happens for almost all $a\pmod{q}$?}
N.M. Korobov \cite{Korobo} showed that for prime $q$  there exists  $a$,
$(a,q)=1$ such that
$$
\max_\nu b_\nu(a)\ll\log q.
$$
Such a result is true for composite $q$ also.
Moreover  Rukavishnikova \cite{ruka1} proved
\begin{theorem}\label{ruuka}
$$
\frac{1}{\varphi (q)}\#\left\{ a\in \mathbb{Z}:\,\,
1\le a \le q,\,\, (a,q) = 1,\,\, \max_{1\le j\le s(a)} b_j (a) \ge T\right\}
\ll \frac{\log q}{T}.
$$
\end{theorem}

Here we 
would like to note that the main results
of Rukavishnikova's papers 
 \cite{ruka1,ruka2} deal with the typical values of the sum of partial quotients  of fractions with a given denominator:
she proves an analog of the law of large numbers.

\subsection{ Exploring folding lemma}

Niederreiter \cite{10} 
proved that  Zaremba's conjecture is true for $q = 2^\alpha,3^\alpha, \,\, \alpha\in \mathbb{Z}_+$ with 
 $\hbox{\got k}=4$, and for $q=5^\alpha$ with $\hbox{\got k}=5$.
His main argument was as follows. If the conjecture is true for $q$ then it is true for $Bq^2$ with bounded integer $B$.
The construction is very simple.
Consider continued fraction
(\ref{exe}) with $b_0(a) = 0$ and its denominator written as a continuant:
$$ q =\langle b_1,b_2,\dots , b_s\rangle
.$$
Define $a^*$ by
$$
 aa^*  \equiv\pm 1\pmod{q}
$$
(the sign $\pm$ should be chosen here with respect to the parity of $s$).
Then
$$
 \frac{a^*}{q}=[0;b_s(a),\dots,b_2(a),b_1(a)]
= \frac{\langle b_1,\dots ,b_{s-1}\rangle}{\langle b_1,\dots , b_s\rangle}
$$
and
$$  q = \langle b_1,\dots ,b_{s-1}\rangle =\langle c_1,c_2,\dots , c_s\rangle,\,\,\, c_j = b_{s-j}.
$$
At the same time if $c_1\ge 2$ then
$$
\frac{q-a^*}{q}=[0;1,c_1-1,\dots,c_s]
=
 \frac{\langle c_1-1,c_2,\dots ,c_{s}\rangle}{\langle 1, c_1-1,c_2,\dots , c_s\rangle}
$$
So we see that
$$
\langle b_1,\dots , b_{s-1}, b_s, X, 1,c_1-1,c_2,\dots , c_s\rangle
=
$$
$$=
 \langle b_1,\dots , b_{s-1}b_s\rangle \langle X,1, c_1-1,c_2,\dots , c_s\rangle
+
\langle b_1,\dots , b_{s-1}\rangle \langle 1,c_1-1,c_2\dots , c_s\rangle=
$$
$$
=
\langle b_1.\dots , b_s\rangle
\langle 1, c_1-1,c_2,\dots , c_s\rangle
\left(X+ \frac{\langle c_1-1,c_2,\dots , c_s\rangle}{\langle 1, c_1-1,c_2,\dots , c_s\rangle
}+
\frac{\langle b_1,\dots , b_{s-1}\rangle}{\langle b_1,\dots , b_s\rangle
}\right)=
$$
$$
= \langle b_1.\dots , b_s\rangle
\langle 1, c_1-1,c_2\dots , c_s\rangle
\left(X+ 1\right)
.$$
This procedure is known as folding lemma. 
                                            
By means of folding lemma Yodphotong and  Laohakosol  showed \cite{yl}
that Zaremba's conjecture is true for $q=6$ and $\hbox{\got k}=6$.
Komatsu \cite{jap2} proved that
Zaremba's conjecture is true for $q=7^{r2^r}, r =1,3,5,7,9,11$ and $\hbox{\got k}=4$.
Kan and Krotkova \cite{kak}
obtained different lower bounds for the  number
$$
f=
\#
\{ a\pmod{p^m}:\,\, a/p^m =[0;b_1,\dots. b_s],\,\, b_j \le p^n\}
$$
of fractions with bounded partial quotients and the denominator of the form $p^n$. 
In particular they proved a bound of the form
$$ 
f \ge C(n)  m^\lambda,\,\,\, C(n), \lambda >0.
$$

 Another applications of folding lemma one can find for example in \cite{bufolding,cohn,akophd,ako}
and in the papers refered there. 
I think that A.N. Korobov proved Niederreiter's result 
concernind powers of 2 and 3 independently in his PhD thesis \cite{akophd}.

\subsection{A result by J. Bourgain and A. Kontorovich (2011)}

Recently J. Bourgain and A. Kontorovich \cite{burr,burra}  
achieved essential progress in Zaremba's conjecture.
Consider the set
$$
{\cal Z}_k (N) :=\{
q\le N:\,\, \exists a \,\,\,\text{such that}\,\,\, (a,q) =1,\,\, \,\, a/q=[0;b_1,...,b_s],\,\,\, b_j \le k\}
$$
( so Zaremba's conjecture means that ${\cal Z}_k (N)= \{1,2,...,N\}$).
In a wonderful paper \cite{burr} they proved 

\begin{theorem}\label{one}
For  $k$ large enough there exists positive $c=c(k)$ such that for $N$ lagre enough
one has
$$
\# {\cal Z }_k(N  )
= N- O(N^{1-c/\log\log N}).
$$
 \end{theorem}
For example it follows from Theorem \ref{one}  that for $k$ large enough the set
$\cup_n {\cal Z}_k (N)$ contains infinitely many prime numbers.

Another result from \cite{burr} is as follows.

\begin{theorem}\label{two}
For  $k =50$ the set
$
\cup_N {\cal Z }_{50}(N  )
$
has positive proportion in $\mathbb{Z}_+$, that is
$$
\#{\cal Z}_{50} (N) \gg N.
$$
 \end{theorem}

These wondeful results follow from right order  upper bound for the integral
\begin{equation}\label{ira}
I_N =
\int_0^1 |S_N(\theta)|^2 d \theta,
\end{equation}
where
$$
S_N(\theta) =\sum_{q\le N} N (q) e^{2\pi i q\theta}, 
$$
and $
N(q) = N_k (q)  $ is the number of integers $a, (a,q)=1$ such that all the partial quotients in the continued fraction
expansion for $\frac{a}{q}$ are bounded by $k$
 
For example positive proportion result follows from the bound
$$
I_N\ll
\frac{S_N(0)^2}{N}
$$
by the Cauchy-Schwarz inequality.

The procedure of estimating of the integral
comes from Vinorgadov's method on estimating of exponential sums with polynomials (Weyl sums).
The main ingredient of the proof (Lemma 7.1 from \cite{burr}) needs 
spectral theory of automorphic forms and follow from a result by
 Bourgain,  Kontorovich and Sarnak  from \cite{gafa}.   

I think that it is possible to simplify the proof given by 
Bourgain and Kontorovich and to avoid the application of a 
difficult result from \cite{gafa}.
Probably for a certain positive proportion result A. Weil's estimates on Kloosterman sums should be enough.\footnote{It was actually done recently by Igor Kan and Dmitrii Frolenkov
\cite{wwk,wwk1,wwk2}. Moreover Kan and Frolenkov improved on some results by Bourgain and Kontorovich on Zaremba conjecture.}

\subsection{ Real numbers with bounded partial quotients}
In this subsection we  formulate some  well-known results concerning real numbers with bounded partial quotients.
We deal with Cantor type sets
$$
F_k =\{ \alpha \in [0,1]: \alpha =[0;b_1,b_2,\dots],\,\, b_j \le k\}
.$$
For the Hausdorff dimension ${\rm dim} \, F_k$ Hensley \cite{HDHEN} proved
$$
r_k = 
{\rm dim}\, F_k =
1- \frac{6}{\pi^2}\frac{1}{k} -
\frac{72}{\pi^4}\frac{\log k}{k^2} + O\left(\frac{1}{k^2}\right),
\,\,\, k \to \infty
$$
Exlicit estimates for ${\rm dim} \, F_k$  for certain values of $k$ one can find in\cite{jeeee}.
Another result by Hensley \cite{HEN,HEN1} is as follows.
For the sums of the values
$$
N_{k}(q) =
\#
\{ a\pmod{q}:\,\, a/q =[0;b_1,\dots. b_s],\,\, b_j \le k\}
$$
considered in the previous subsection
Hensley proved
\begin{equation}\label{he}
\sum_{q\le Q} N_k(q) \sim {\rm constant} \times Q^{2r_k}
,\,\,\, Q \to +\infty.
\end{equation}

We need  a
corollary to (\ref{he}).
Consider the set 
$$
B(k,T)=
\{ a\pmod{q}:\,\, a/q =[0;b_1,\dots. b_\nu,\dots , b_s],
\,\,
\langle b_1,...,b_\nu \rangle \le T \,\,\Longrightarrow\,\,
\,\, b_j \le k\}.
$$
Then for $T\ll \sqrt{q}$ one has
$$
\#
B(k,T)\asymp_k
q^{r_k}.
$$

\subsection{Zaremba's conjecture and points on modular hyperbola}

 When we are speaking about ``modular hyperbola`` we are interested in the distribution of the points from the set
$$
\{ (x_1,x_2) \in \mathbb{Z}_q^2:\,\,\,
x_1x_2 \equiv \lambda \pmod{q}\}
.$$
For a wonderful survey we would like to refer to Shparlinski \cite{spa}.

{\bf Proposition.} \,\,\,{\it
 Suppose that for $
T \ge C \sqrt{q}$ there exist $x_1,x_2 \pmod{q}$ such that
$$
x_1,x_2 \in B(k,T),
$$
and
$$
x_1x_2 \equiv 1\pmod{q}.
$$
Then  Zaremba's conjecture is is true
with a certain $\hbox{\got k}$ depending on $k$ and $C$.
}

From A. Weil's bound for complete Kloosterman sums
we know that the points on modular hypebola are
uniformily distributed in boxes of the form
$$
I_1\times I_2,\,\,\, I_\nu = [X_\nu, X_\nu + Y_\nu]
$$
provided
$$
Y_1\times Y_2 \ge q^{3/2+\varepsilon}.
$$
So we have the following

{\bf
Corollary.}\,\,\,{\it
Suppose that
\begin{equation}\label{bett}
\beta_1+\beta_2 <\frac{1}{4}.
\end{equation}
Then there exist $x_1,x_2 \pmod{q}$ such that
$$
x_1x_2 \equiv 1\pmod{q},\,\,\,\,\,\,
x_1 \in B(k,T_1),
\,\,\, x_2 \in B(k,T_2)
$$
and
$$
T_1\asymp q^{\beta_1},\,\,\, T_2 \asymp q^{\beta_2},$$}

 By means of application of bounds for incomplete Kloosterman sums Moshchevitin \cite{apb}
proved the following result.

\begin{theorem}\label{mot}
Put
$$
\omega_1 = \omega_1 (\beta_1 ) = \frac{(4r_k\beta_1 -1)(1 -
2\beta_1 )}{8}
 .
$$
 Let $ q=p$ be prime and
$$
\frac{1}{4r_k} < \beta_1 <\frac{1}{2},\,\,\,  0 < \beta_2
<\omega_1 .
$$
Put $T_j = p^{\beta_j},j = 1,2$.
Then there exist $$
x_1 \in B(k,T_1),
\,\,\, x_2 \in B(k,T_2)
$$
such that
$$
x_1x_2\equiv 1\pmod{p}.
$$
\end{theorem}

Theorem \ref{mot}
improves on the Corollary above in the case of prime $q=p$, as we can take $\beta_,\beta_2$ in the range
$$
\frac{1}{2} -\varepsilon <\beta_1+\beta_2 <\frac{1}{2}
,$$ 
 which is not considered in (\ref{bett}).
However  in Theorem \ref{mot} there is no symmetry between $\beta_1$ and $\beta_2$.
It gives no good result for $\beta_1 = \beta_2$.

\subsection{Numbers with missing digits } 

In this section we will show that
if one considers instead of ''non-linear'' fractal-like sets $B(k,T)$ a more simple fractal-like set,
the corresponding problem becomes much more easier.

For positive integers
 $s, k$ we consider sets
$$D = \{d_0,..., d_k \}
,\,\, 0 =d_0 < d_1 < ... < d_k
<s,\hspace{2mm} 1\le k \le s-2,\,\,\,(d_1,...,d_k) = 1, $$
$$
 K_s^D (N) = \{ x\in \mathbb{Z}_+:
\,\,\, x < N,
  \hspace{2mm} x= \sum_{j=0}^h \delta
_j s^j, \hspace{2mm} \delta _j \in D \} \label{K}
.$$

We are interested in properties of elements of $K_s^D(N)$ modulo $q$.

 In \cite{mmm}  by means of A. Weil's bounds for Kloosterman sums the following result was proven.
Let $p$ be prime.
 Under certain natural conditions on $s,D$
for any $\lambda\pmod{p}$
 there exist $x_1, x_2 \in  K_s^D (p) $
such that 
\begin{equation}\label{twaila}
x_1x_2 \equiv \lambda\pmod{p}.
\end{equation}
One can compare this result with proposition from Subsection 5.5.

We conclude this subsection by mentioning an interesting open  problem formulated by Konyagin \cite{ook}.
The question is as follows. 

{\it
 Suppose that 
$(s,q)=(d_0,...,d_k)=1, k \ge 2.$
Is it true that for some large $\sigma>0$
for
\begin{equation}\label{jool}
N \ge q^\sigma\end{equation}
 for any $\lambda \pmod{q}$ there exists 
$
x \in  K_s^D (N)
$
 such that
$
x\equiv \lambda \pmod{q}\,?
$
}

Konyagin \cite{ook} showed that the answer is ``yes``  for {\it almost all} $q$ (a simple variant of large sieve argument).
In the same paper he showed that the conclusion is true  if we replace
the condition (\ref{jool}) by $ N \ge \exp(\sigma \log q \log\log q)$.
Some related topics were considered by the author in \cite{x1,x2,x3}.
In \cite{mmm} it is shown that under the condition (\ref{jool}) with $\sigma$ large enough
and $q=p$ prime  for any $\lambda$ there exist $x_1, x_2 \in 
K_s^D (N)$ satisfying (\ref{twaila}).

\subsection{Discrepancy bounds}

For $(a,q) =1$ consider the discrepansy $D(a,q)$ of the finite sequence of points
\begin{equation}\label{see}
\xi_k=
\left(\frac{k}{q},\left\{\frac{ak}{q}\right\}\right),\,\,\, 0\le k \le q-1
.
\end{equation}
It is defined as
$$
D(a,q) =
\sup_{\gamma_1,\gamma_2\in (0,1)}
\left|
\#\{k:\xi_k \in [0,\gamma_1)\times [0,\gamma_2)\} - q\gamma_1\gamma_2\right|.
$$
Here we do not want to discuss the foundations an major results of the theory of uniformly distributed sequences;
we refer to books \cite{nied} and \cite{dete}. 

It is a well-known fact that  
\begin{equation}\label{summa}
D(a,q) \ll
\sum_{j=1}^{s(a)} b_j (a)
\end{equation}
where $b_j$ are partial quotients from (\ref{exe}).

If Zaremba's conjecture is  true then for any $q$ there exists $a$ coprime to $q$
such that
\begin{equation}\label{conti}
D(a,q)\ll \log q.
\end{equation}
  
Larcher \cite{Lar} proved that  for any $q$ there exists $a$ coprime to $q$ such that
\begin{equation}\label{lara}
D(a,q) \ll \frac{q}{\varphi (q)} \,\log q \log\log q.
\end{equation}
This bound is optimal up to the factor $\log\log q$.
In fact from Rukavishnikova's results \cite{ruka1,ruka2} we see that (\ref{lara})
holds for almost all $a$ coprime to $q$.

However Zaremba's conjecture is still open, and we do not know if  for a given  $q$ one can get (\ref{conti})
instead of (\ref{lara}), for some $a$.

Ushanov and Moshchevitin \cite{Um}
proved the following result.

\begin{theorem}\label{repeter}
{ Let $p$ be prime, $U$ be a multiplicative  subgroup in $\mathbb{Z}_p^*.$  
For  $v \neq 0$ we consider the set $R = v\cdot U$ and let
\begin{equation}\label{burrt}
\#R \geq 10^8 p^{7/8}\log^{5/2}{p}.
\end{equation}	
Then there exists an element $a \in R,$ $a/p=[b_1, b_2, \cdots, b_l],$ $b_i = b_i(a),$ $l = l(a)$ with
$$
\sum_{i=1}^{l} b_i \leq 500\log{p}\log\log{p},
$$
and hence
$$
D(a,p) \ll \log p\log\log p.
$$
}
\end{theorem}
The proof uses Burgess' inequality for character sums.

An open problem here is as follows. {\it 
Is it possible to replace exponent $7/8$ in the condition (\ref{burrt}) by a smaller one}?

Recently Professor Shparlinski informed me that 
Chang \cite{changa} essentially repeated the  result of
Theorem \ref{repeter}.
Moreover in \cite{changa} an analog of Rukavishnikova's Theorem \ref{ruuka} is proved for multiplicative subgroups
$\pmod{p}$  of cardinality $\gg p^{7/8+\varepsilon}$.

A multidimensional version of Theorem \ref{repeter} was obtained by
Ushanov \cite{ushaa}. It is related to Theorem \ref{biko} below.

\subsection{Discrepancy bounds: multidimensional case}

We would like to conclude this section by mentioning a wonderful recent result due to Bykovskii \cite{Byk,Byk1} dealing with multidimensional analog of the sequence 
(\ref{see}).
We consider positive integer $s$ and integers $q\ge 1;a_1=1,a_1,...,a_s$.
The study of the distribution of the sequence
$$
\xi_k = \left(
\frac{a_1k}{q},\left\{\frac{a_2k}{q}\right\},...,\left\{\frac{a_sk}{q}\right\}\right),\,\,\, 0\le k \le q-1
$$
started with the works of Korobov \cite{korodan} and Hlawka \cite{hlaa}.
We are interested in upper bounds for the discrepancy
$$
D(a_1,...,a_s;q) =
\sup_{\gamma_1,...,\gamma_s\in (0,1)}
\left|
\#\{k:\xi_k \in [0,\gamma_1)\times \cdots\times [0,\gamma_s)\} - q\gamma_1\cdots\gamma_s\right|.
$$
The  upper bound
$$
\min_{(a_1,...,a_s) \in \mathbb{Z}^s}
D(a_1,...,a_s;q)\ll_s
{(\log q)^s}
$$
was proved by Korobov \cite{Korobo,koroboumn} for prime $q$ and by Niederreiter \cite{optima}
for composite $q$.

Bykovskii \cite{Byk,Byk1} proved the following result.

\begin{theorem}\label{biko}
For $s \ge 2$ one has  and for any positive 
integer $q$ one has
\begin{equation}\label{byyk}
\min_{(a_1,...,a_s) \in \mathbb{Z}^s}
D(a_1,...,a_s;q)\ll_s
{(\log q)^{s-1} \log\log q}
.\end{equation}
\end{theorem}
For $s=2$  and prime $q$ this result coincides with Larcher's inequality (\ref{lara}) discussed in the previous subsection.
However Larcher's proof is based on the inequality 
(\ref{summa})
while Bykovskii's proof is related to analytic argument (this argument is quite similar for the case $s=2$ and
the general case $ s\ge 2$) and to consideration of relative minima of lattices.
It happened that the proof of Bykovskii's result is not extremely difficult.
It is related to a paper by Skriganov \cite{skri} and the previous paper by Bykovskii \cite{byda}.
The method developed by Bykovskii may find applications in other problems (see for example \cite{frol}).

A famous well-known conjecture is that
$$
\min_{(a_1,...,a_s) \in \mathbb{Z}^s}
D(a_1,...,a_s;q)\ll_s
{(\log q)^{s-1} }
,$$
that is that the factor $\log\log q$ in (\ref{byyk}) may be thrown away. This conjecture seems to be very difficult.

\section{Minkowski question mark function}

\subsection{Definition of Minkowski function}
If real  $x=[0;a_1,...,a_t,...] \in [0,1] $ is represented as a regular continued fraction with natural partial quotients, then
Minkowski question mark function $?(x)$ is defined as follows:
$$
?(x) = \frac{1}{2^{a_1-1}} - \frac{1}{2^{a_1+a_2-1}}+ ...+ \frac{(-1)^{n+1}}{2^{a_1+...+a_n-1}}+ ...  
$$
(in the case of rational $x$ this sum is finite).
It is a well known fact that $?(x)$ is a continuous strictly increasing function. By the Lebesgue theorem
it has finite derivative almost everywhere in $[0,1]$.
Moreower  the derivative  $?'(x)$, if exists (in finite or infine sense) can have only two values -  $0$  or $+\infty$.
For more results  we refer to papers \cite{alk,alk1,den,duma,MI,KI,PARA,paradizo,Sal}.

It will be important for us to recall the definition of Stern-Brocot sequences  $ F_n $, $ n=0, 1, 2, \dots $. 
For $n=0$ one has
  $$ F_0=\{0, 1\}=\left\{\frac{0}{1}, \frac{1}{1}\right\} .$$
Suppose that the sequence $ F_n $ is written in the increasing order
 $$0=\xi_{0, n}<\xi_{1, n}< \dots <\xi_{N \left( n
  \right), n}=1, N(n)=2^{n},
\,\,\,\, \xi_{j, n} =\frac{p_{j, n}}{q_{j, n}},\,\,\, ( p_{j, n},q_{j, n}) = 1.$$
Then the sequence  $F_{n+1}$ is defined as   $$ F_{n+1} = F_n \cup
Q_{n+1} $$ where 
$$ Q_{n+1}=\left\{
\frac{p_{j, n}+ p_{j+1, n} }{q_{j, n}+ q_{j+1, n}},\,\,\,
 j=0, \dots , N(n)-1\right\}.  $$  
Note that for the number of elements in $F_n$ one has
$$
\#F_n = 2^{n}+1.
$$
The Minkowski question mark function $?(x)$ is the limit distribution function for the Stern-Brocot sequences:
$$
?(x)  =\lim_{n\to \infty} \frac{\# \{ \xi\in F_n :\,\,\,\xi \le x\}}{2^n +1}.
$$

\subsection{Fourier-Stieltjes coefficients}

Here we would like to mention a famous open problem  by R. Salem \cite{Sal}: {\it
to prove or to disprove that for Fourier-Stieltjes
coefficients $d_n, n \in \mathbb{N}$  of $?(x)$ one has}
$$
d_n = \int_0^1 \cos (2\pi n x) \, {\rm d}?(x) \to 0,\,\,\, n  \to \infty.
$$
Certain results related to this problem were obtained by G. Alkauskas \cite{alk,alk1}.

\subsection{Two simple questions}
Here we formulate two open questions.

1. One can see that
$$
?(0)= 0,\,\,\,\,\,\, ?\left(\frac{1}{2}\right) = \frac{1}{2},\,\,\,\,\,\,
?(1) = 1.
$$
Moreover
$$
?'(0) = ?'\left(\frac{1}{2}\right)
=?'(1) =0,
$$
as in any rational point the question mark function has zero derivative.
By continuouity agrument we see that there exist two points
$$
x_1 \in \left(0,\frac{1}{2}\right),\,\,\,\,\,
x_2 \in \left(\frac{1}{2}, 1\right)
$$
such that
$$
?(x_i) = x_i,\,\,\, i=1,2.
$$
So we see that the equation
\begin{equation}\label{woo}
?(x) = x,\,\,\,\, x \in [0,1]
\end{equation}
has at least five solutions.
The question is  if equation (\ref{woo}) has
{\it exactly} five solutions.

2. Consider the function $m (x) $ inverse to $?(x)$. As
$?(\xi_{j,n}) = \frac{j}{2^n}$
we see that $m\left(\frac{j}{2^n}\right) =\xi_{j,n}$.\
Then by the Koksma inequality (see \cite{nied}) we have
$$\left|
\frac{1}{2^n}\,
\sum_{j=1}^{2^n} \left(\xi_{j,n} - \frac{j}{2^n}\right)^2 -
\int_0^1 ( m(x) - x)^2 {\rm d}x \right|
\le \frac{D_n \cdot V}{2^n},
$$
where $D_n$ is the discrepancy of the sequence
$$
\frac{j}{2^n},\,\,\, 1\le j \le 2^n,
$$
$ 0<D_n \le 1$ 
and $V\le 4$ is the variation of the function $x\mapsto (m(x)-x)^2$.
One can easily see that 
\begin{equation}\label{iii}
\int_0^1 ( m(x) - x)^2 {\rm d}x =
\int_0^1 ( ?(x) - x)^2 {\rm d}x >0
.
\end{equation}
That is why we have
\begin{equation}\label{fraa}
 \sum_{j=1}^{2^n} \left(\xi_{j,n} - \frac{j}{2^n}\right)^2 =
2^n \int_0^1 ( ?(x) - x)^2 {\rm d}x +R_n,\,\,\,
|R_n|\le 4.
\end{equation}
The question is as follows. {\it Is it true that for the remainder in (\ref{fraa}) one has
$R_n \to 0 $ as $ n \to \infty$}?

This question is motivated by the famous Franel theorem. Instead of Stern-Brocot sequence $
F_n$ consider Farey series ${\cal F}_Q$ which consist of all rational numbers $p/q\in [0,1], (p,q) = 1$
with denominators $\le Q$. Suppose that ${\cal F}_Q$ form an increasing sequence
$$
1= r_{o,Q}<r_{1,Q}<...< r_{j,Q}<r_{j+1,Q} <...< r_{\Phi (Q),Q}= 1, \,\,\,\,
\Phi (Q) 
=\sum_{q\le Q}\varphi (q)
$$
(here $\varphi (\cdot )$ is the Euler totient function).
Then for the limit distribution function one has
$$
\lim_{Q\to \infty} \frac{\# \{ r\in {\cal F}_Q :\,\,\, r \le x\}}{\Phi (Q)+1} =x
$$
and the integral similar to (\ref{iii}) is equal to zero.
Franel's theorem (see \cite{Landau}) states that the asymptotic formula
$$
 \sum_{j=1}^{\Phi (Q)} \left(r_{j,Q} - \frac{j}{\Phi (Q)}\right)^2 = O_\varepsilon (Q^{-1+\varepsilon}) 
, \,\,\, Q\to \infty.
$$
for all positive $\varepsilon$ 
is equivalent to Riemann Hypothesis.
In fact the well-known asymptotic equality
$$
\sum_{n\le Q} \mu (n) = o(Q), \,\,\,
Q\to \infty
$$
leads to
$$
\sum_{j=1}^{\Phi (Q)} \left(r_{j,Q} - \frac{j}{\Phi (Q)}\right)^2 = o(1) 
, \,\,\, Q\to \infty.
$$
A analogous formula for $R_n$ from (\ref{fraa}) is unknown, probably.

\subsection{Values of derivative}

It is a well-known fact that if for $x \in [0,1]$ the derivative $?'(x)$ exists then
$?'(x) = 0$ or $?'(x) = +\infty$.

For a real irrational $x$ represented as a condinued fraction expansion
$x=[a_0;a_1,...,a_n,..]$
we consider the sum of its first partial quotients
$S_x(t) =a_1+...+a_t$.
Define
$$
\kappa_1 
=
\frac{2 \log \frac{1+\sqrt{5}}{2}}{\log 2} =  1.388^+,\,\,\,\,
\kappa_2=
\frac{4L_5-5L_4}{L_5 -L_4}= 4.401^+,\,\,\,\,
L_j
= \log \frac{j+\sqrt{j^2+4}}{2} - j\frac{\log 2}{2}.
$$
Improving on  results by 
 Paradis,  Viader and Bibiloni \cite{paradizo}, Kan, Dushistova and Moshchevitin \cite{mdk}
proved the following four theorems.

\begin{theorem}\label{01}

(i) Assume for an irrational number $x$ there exists such a constant $C$
that for all natural $t$ 
one has
$$
S_x(t)\le \kappa_1t+\frac{\log t}{\log 2}+C. 
$$
Then $?'(x)$ exists and $?'(x)=+\infty$.

(ii) Let $\psi (t)$ be an increasing function such that $\lim_{t\to +\infty} \psi (t) =+\infty$. Then there exists such an irrational number $x\in
(0,1)$ that $?'(x)$ does not exist and for any $t$ one has
$$
S_x(t)\le \kappa_1t+\frac{\log t}{\log 2} +\psi(t). 
$$
\end{theorem}

\begin{theorem}
\label{02}

Let for an irrational number $x\in (0,1)$ the derivative
$?'(x)$ exists and $?'(x)=0$. Then
for any real function $\psi =\psi (t)$ under conditions
$$
\psi (t) \ge 0,\,\,\,\, \psi (t) = o\left(\frac{\log\log t}{\log t}\right),\,\, t\to\infty
$$
there exists  $T$ depending on $\psi$ such that for all $t\ge T$ one has
 $$
\max_{u\le t}\left(S_x(u)-\kappa_1 u\right) \ge \frac{ \sqrt{2\log \lambda_1-\log 2} }{\log 2}\cdot\sqrt{t\log
t}\cdot (1-t^{-\psi (t)}).
$$

(ii) There exists such an irrational  $x\in (0,1)$ that $?'(x)=0$ and for all large enough $t$ 
$$
S_x(t)-\kappa_1 t \le 
 \frac{ \sqrt{16\log \lambda_1-8\log 2} }{\log 2}\cdot\sqrt{t\log
t}\cdot \left(1+2^5\left(\frac{\log\log t}{\log t}\right)\right).
$$
\end{theorem}

\begin{theorem}\label{03}

(i) Assume for an irrational number $x$ there exists such a constant $C
$ that for all natural $t$   
one has
$$
S_x(t)\ge \kappa_2t-C. $$
Then $?'(x)$ exists and $?'(x)=0$.

(ii) Let $\psi (t)$ be an increasing function such that $\lim_{t\to +\infty} \psi (t) =+\infty$. Then there exists such an irrational number $x\in
(0,1)$ that $?'(x)$ does not exist and for any $t$ 
we have
$$
S_x(t)\ge \kappa_2t-\psi(t). 
$$
\end{theorem}

\begin{theorem}\label{04}

(i) Assume for an irrational number $x\in (0,1)$ the derivative $?'(x)$ exists and 
$?'(x)=+\infty$. Then for any large enough $t$  
one has
$$
\max_{u\le t}\left( \kappa_2 u -S_x(u)\right) \ge  \frac{ \sqrt{t}}{10^8}.
$$

(ii) There exists such an irrational $x\in (0,1)$ that $?'(x)=+\infty$ and for large enough $t$ 
we have
$$
\kappa_2 t - S_x(t)\le 200 \sqrt{t}.
$$
\end{theorem}

A weaker result is due to Dushistova and Moshchevitin \cite{duma}.
All these results are related to deep analysis of sets of values of continuants.
Theorems \ref{01} and \ref{03} are optimal, but Theorms \ref{02} and \ref{04}
are not optimal. It should be interesting to prove optimal bounds  related to this two last theorems.

Here we should note that the paper \cite{mdk} contains some other results related to sets of real numbers 
with bounded partial quotients.
Probably some of these results may be improved.

It is possible to prove that for any $\lambda$ from the interval
 $$
 \kappa_1\le \lambda\le \kappa_2
 $$
 there exist irrationals $x,y,z \in [0,1]$ such that
 $$
\lim_{t\to \infty}\frac{S_x(t)}{t}=\lim_{t\to \infty}\frac{S_y(t)}{t}=\lim_{t\to \infty}\frac{S_z(t)}{t}=\lambda
$$
and $?'(x) =0, ?'(y) = +\infty$, but $?'(z)$ does not exist.

Theorems \ref{01} - \ref{04}
show that
\begin{equation}\label{w1}
\kappa_1=\sup
\left\{
\kappa\in \mathbb{R}:\,\,
\limsup_{t\to \infty}\frac{S_t(x)}{t} <\kappa\,\,\,
\Longrightarrow\,\,\, ?'(x) =+\infty \right\},
\end{equation}
\begin{equation}\label{w2}
\kappa_2=\inf
\left\{
\kappa\in \mathbb{R}:\,\,
\liminf_{t\to \infty}\frac{S_t(x)}{t} >\kappa\,\,\,
\Longrightarrow\,\,\, ?'(x) = 0 \right\}.
\end{equation}

\subsection{Denjoy-Tichy-Uitz family of functions}

There are various generalizations of the Minkowski question mark function $?(x)$. One of them was considered by Denjoy \cite{den}
and rediscovered  by Tichy and Uitz \cite{TU}.

For $\lambda \in (0,1)$ we define for $x\in [0,1]$ a function $g_\lambda (x)$ in the following way.
Put
$$g_\lambda (0) =0,\,\,\,\,  \,\,\,\,
g_\lambda(1) = 1
.$$
Then if $g_\lambda $ is defined for two neighbooring  Farey fractions $\frac{a}{b}< \frac{c}{d}$,  we put
$$
g_\lambda \left(
\frac{a+c}{b+d}\right) =
(1-\lambda) g_\lambda \left(\frac{a}{b}
\right)+
\lambda
g_\lambda \left(\frac{c}{d}\right).
$$
So we define $g_\lambda (x)$ for all rational $x \in [0,1]$.
For irrational $x$ we define $g_\lambda (x)$ by continuouty.

The family $\{g_\lambda \}$ constructed consist of sungular functions.
One can easily see that $g_{1/2} (x) = ?(x)$, that is the Minkowski question mark function is a member of this family.
Hence $g_{1/2}$ has a clear arithmetic nature: it is the limit distribution function for Stern-Brocot sequences $F_n$.
Recall that
$$
F_n =\{ x=[0;a_1,...,a_t]:\,\, a_1+...+a_t = n+1\}.
$$

Zhabitskaya \cite{z} find out that
for $\lambda = \frac{3-\sqrt{5}}{2}$
the function $g_\lambda (x)$ is the limit distribution function for the sequences 
$$
\Xi_n = \{ x=[[1;b_1,...,b_l]]:\,\,\, b_1+...+b_l = n+1\}
$$
associated with ``semiregular'' continued fractions
$$
[[1;b_1,b_2,b_3,...,b_l]]=
1-
\frac{1}{\displaystyle{b_1-\frac{1}{\displaystyle{b_2-
\frac{1}{\displaystyle{b_3-\dots -
\displaystyle{\frac{1}{b_{l}}} }}}}}} , \,\,\,\,\, b_j \in \mathbb{Z},\,\,\, b_j \ge 2
.$$

It happens that disrtibution functions of some other sequences associated with special continued fractions
 do not belong to the family
$\{g_\lambda \}$  (see \cite{z1,z2}).

It is interesting to find other values of $\lambda$ for which the function  $g_\lambda(x)$ is associated with an
explicit and natural object.

Another open problem related to the family 
$\{g_\lambda \}$  is as follows. Analogously to (\ref{w1},\ref{w2}) we define
$$
\kappa_1 (\lambda)=\sup
\left\{
\kappa\in \mathbb{R}:\,\,
\limsup_{t\to \infty}\frac{S_t(x)}{t} <\kappa\,\,\,
\Longrightarrow\,\,\, g_\lambda'(x) =+\infty \right\},
$$
$$
\kappa_2(\lambda)=\inf
\left\{
\kappa\in \mathbb{R}:\,\,
\liminf_{t\to \infty}\frac{S_t(x)}{t} >\kappa\,\,\,
\Longrightarrow\,\,\, g_\lambda'(x) = 0 \right\}.
$$
The study of the functions
$$
\kappa_j (\lambda),\,\,\, 0<\lambda <1
$$
has never been made.
Recently D. Gayfulin calculated the values $\kappa_{j} \left(\frac{3-\sqrt{5}}{2}\right),\, j = 1,2$.

\end{document}